\documentclass[a4paper,11pt]{article}

\usepackage{anysize}
\usepackage{amsthm}
\usepackage{amssymb}
\usepackage{amsmath}
\usepackage{fancybox}
\usepackage{amsfonts}
\usepackage[T1]{fontenc}
\usepackage{latexsym}
\usepackage{makeidx}
\usepackage{cite}
\usepackage[numbers]{natbib}

\usepackage[pdftex]{graphicx}
\usepackage[pdftex]{color}

\usepackage[all]{xy}

\newtheorem{teo}{Theorem}
\newtheorem{lem}{Lemma}
\newtheorem{pro}{Proposition}
\newtheorem{cor}{Corollary}
 
\newcommand{\Fd}{\longrightarrow}
\newcommand{\fd}{\rightarrow}

\newcommand{\inc}{\subset}

\newcommand{\iso}{\cong}

\newcommand{\pt}{\forall}

\newcommand{\ba}{\overline}

\newcommand{\al}{\alpha}

\newcommand{\del}{\delta}

\newcommand{\gam}{\gamma}

\newcommand{\Lam}{\Lambda}
\newcommand{\Sig}{\Sigma}
\newcommand{\cit}{\theta}

\newcommand{\Om}{\Omega}

\newcommand{\Z}{\mathbb{Z}}
\newcommand{\N}{\mathbb{N}}

\newcommand{\R}{\mathbb{R}}
\newcommand{\C}{\mathbb{C}}
\newcommand{\K}{\mathbb{K}}

\newcommand{\Sa}{\mathbb{S}}

\newcommand{\pa}{\partial}

\newcommand{\p}{\grave{}}

\newtheorem{remark}{Remark}[section]

\def\pf{\par\noindent {\em Proof.}~\par\noindent}

\def\lim{\mathop{\mbox{\normalfont lim}}\limits}
\def\sd{\mathop{\mbox{\normalfont sdet}}\limits}
\def\str{\mathop{\mbox{\normalfont str}}\nolimits}
\def\Mat{\mathop{\mbox{\normalfont Mat}}\nolimits}
\def\GL{\mathop{\mbox{\normalfont GL}}\nolimits}
\def\Sym{\mathop{\mbox{\normalfont Sym}}\nolimits}
\def\tr{\mathop{\mbox{\normalfont tr}}\nolimits}
\def\Pin{\mathop{\mbox{\normalfont Pin}}\nolimits}
\def\Spin{\mathop{\mbox{\normalfont Spin}}\nolimits}
\def\pf{\par\noindent {\em Proof.}~\par\noindent}

\def\pa{\partial}

\begin{document}

\date{}

\title{The Spin Group in Superspace}
\small{
\author
{Hennie De Schepper, Al\'i Guzm\'an Ad\'an, Frank Sommen}
\vskip 1truecm
\date{\small  Clifford Research Group, Department of Mathematical Analysis, Faculty of Engineering
and Architecture, Ghent University, Krijgslaan 281, 9000 Gent, Belgium. \\
\texttt{hennie.deschepper@ugent.be}, \texttt{ali.guzmanadan@ugent.be}, \texttt{franciscus.sommen@ugent.be}}

\maketitle

\begin{abstract} 
\noindent There are two well-known ways of describing elements of the rotation group SO$(m)$. First, according to the Cartan-Dieudonné theorem, every rotation matrix can be written as an even number of reflections. And second, they can also be expressed as the exponential of some anti-symmetric matrix. \\
{In this paper, we study similar descriptions of a group of rotations SO${}_0$ in the superspace setting. This group can be seen as the action of the functor of points of the orthosymplectic supergroup $\textup{OSp}(m|2n)$ on a Grassmann algebra. }
While still being connected, the group SO${}_0$  is thus no longer compact. As a consequence, it cannot be fully described by just one action of the exponential map on its Lie algebra. Instead, we obtain an Iwasawa-type decomposition for this group in terms of three exponentials acting on three direct summands of the corresponding Lie algebra of supermatrices. \\
At the same time, SO${}_0$ strictly contains the group generated by super-vector reflections. Therefore, its Lie algebra is isomorphic to a certain extension of the algebra of superbivectors.  This means that the Spin group in this setting has to be seen as the group generated by the exponentials of the so-called extended superbivectors in order to cover SO${}_0$. We also study the actions of this Spin group on supervectors and provide a proper subset of it that is a double cover of SO${}_0$. Finally, we show that every fractional Fourier transform in n bosonic dimensions can be seen as an element of this spin group.

\noindent

\vspace{0.3cm}

\small{ }

\noindent
\textbf{}
\end{abstract}

\section{Introduction}
Supermanifolds and in particular superspaces play an important r\^ole in contemporary theoretical physics, e.g.\ in the particle theory of supersymmetry, supergravity or superstring theories, etc. Superspaces are equipped with both a set of commuting co-ordinates and a set of anti-commuting co-ordinates. From the mathematical point of view, they have been studied using algebraic and geometrical methods. Some pioneering references in the development of analysis on superspace are \cite{Berezin:1987:ISA:38130, MR565567, MR778559, MR0580292, MR760837}. For more modern treatments we address the reader to \cite{MR2840967, MR1701597, MR1632008, MR2069561}. More recently, harmonic and Clifford analysis have been extended to superspace by introducing some important differential operators (such as Dirac and Laplace operators) and by studying special functions and orthogonal polynomials related to these operators, see e.g.\ \cite{Bie2007, de2007clifford, MR2344451, 1751-8121-42-24-245204, MR2521367, MR2683546}. 

In classical harmonic and Clifford analysis in $\R^m$, the most important invariance group is the set of rotations SO$(m)$, i.e.\ the connected group of $m\times m$ real matrices leaving invariant the Euclidean inner product $ \langle \underline{x},\underline{y}\rangle=-\frac{1}{2}\left(\underline{x}\underline{y} + \underline{y}\underline{x}\right) = \sum_{j=1}^m x_jy_j $, $\underline{x},\underline{y}\in \R^m$.
 Every rotation in SO$(m)$ can be written as the exponential of some anti-symmetric matrix and vice versa, each one of such exponentials is a rotation in $\R^m$. The Lie algebra of SO$(m)$ is given by the set $\mathfrak{so}(m)$ of anti-symmetric matrices. 
 
On the other hand, the group SO$(m)$ can also be described by means of the spin group $\mbox{Spin}(m):=\left\{\prod_{j=1}^{2k} \underline{w}_j: k\in\N, \underline{w}_j\in \Sa^{m-1}\right\}$,
where $\Sa^{m-1}=\{\underline{w}\in\R^m: \underline{w}^2=-1\}$ denotes the unit sphere in $\R^m$. The relation between $\mbox{Spin}(m)$ and SO$(m)$ is easily seen through the Lie group representation $h: \mbox{Spin}(m)\fd \mbox{\normalfont SO}(m)$
\[
h(s)[\underline{x}]=s\underline{x}\ba{s}, \hspace{1cm} s\in \mbox{Spin}(m), \; \underline{x}\in \R^{m}, 
\]
which describes the action of every element of SO$(m)$ in terms of Clifford multiplication. It esily follows from the above representation that $\mbox{Spin}(m)$ constitutes a double cover of SO$(m)$.
Such a description of the spin group follows from the Cartan-Dieudonn\'e theorem which states that every orthogonal transformation in an $m$-dimensional symmetric bilinear space can be written as the composition of at most $m$ reflections. Basic references for this setting are \cite{MR1169463, MR1130821}.

In this paper we study the similar situation in the framework of Clifford analysis in superspace, where the Cartan-Dieudonn\'e theorem is no longer valid. The main goal is to properly define {a spin group} as a set of elements describing every super-rotation through Clifford multiplication.  To that end, we consider linear actions on supervector variables using both commuting and anti-commuting coefficients in a Grassmann algebra ${\Lam}(f_1,\ldots, f_N)$. This makes it possible to study the group of supermatrices leaving the inner product invariant and to define in a proper way the spin group in this case. It is worth noticing that the superstructures are absorbed by the Grassmann algebras leading to classical Lie groups and Lie algebras instead of supergroups or superalgebras. 

Finite dimensional Grassmann algebras are the most simple coefficient algebras to consider if one wants to study linear actions mixing bosonic and fermionic variables. Moreover, the use of this coefficient algebra is already rich enough for the purposes of Clifford analysis. Indeed, in \cite{DS_Guz_Somm4} it was shown that this group theoretical approach using Grassmann coefficients offers a new way of describing the $\mathfrak{osp}(m|2n)$-invariance of the super Dirac operator (see also \cite{MR3375856}). This approach also leads to new insights in the study of linear actions (supermatrices) on supervectors. For example, the rotation group is no longer compact and it strictly contains the subgroup generated by an even number of supervector reflections (see Sections \ref{InvIPS} and \ref{SecSpGSP}).


In a more algebraic geometric approach, one can generalize these groups by using the language categories. For example, one can define the group SO${}_0$ of superrotations via the functor of $A$-points that maps the category $\R$-{\bf Salg} of real commutative superalgebras to the category {\bf Grp} of groups as follows
\begin{align*}
\R\mbox{-{\bf Salg}} \fd \mbox{\bf Grp}& \, :  \;A \mapsto \mbox{SO}_0(A), & A&\in \R\mbox{-{\bf Salg}}.
\end{align*}
Here,  $\mbox{SO}_0(A)$ denotes the group of $(m|2n)$-supermatrices satisfying the defining conditions of  $\mbox{SO}_0$ (see Section \ref{SupRot-Grp}) but with entries in the superalgebra $A$. This functor can be proved to be representable since superrotations are solutions of some set (super) polynomial equations (see e.g. \cite[Chpt.\ 11]{MR2840967} and \cite{MR2555977} for other examples of representations of affine supergroups). {In fact, this functor is isomorphic to the functor of points of the orthosymplectic supergroup OSp$(m|2n)$.} However, the representability of the functor associated to the Spin group defined in this paper (see Section \ref{SecSpGSP}) seems to be more difficult to prove. We shall address this problem in future work.

The paper is organized as follows. We start with some preliminaries on Grassmann algebras, Grassmann envelopes and supermatrices in Section \ref{Prem}. In particular, we carefully recall the notion of an exponential map for Grassmann numbers and supermatrices as elements of finite dimensional associative algebras. In Section \ref{SsFwGC} we briefly describe the Clifford setting in superspace leading to the introduction of the Lie algebra of superbivectors. An extension of this algebra plays an important r\^ole in the description of the spin group in this setting. The use of the exponential map in such an extension 
necessitates the introduction of the corresponding tensor algebra. Section \ref{InvIPS} is devoted to the study of the invariance of the ``inner product'' of supervectors. There, we study several groups of supermatrices and in particular, the group of superrotations SO${}_0$ and its Lie algebra $\mathfrak{so}_0$, which combine both orthogonal and symplectic structures. We prove that every superrotation can be uniquely decomposed as the product of three exponentials acting in some special subspaces of $\mathfrak{so}_0$. Finally, in section \ref{SecSpGSP} we study the problem of defining the spin group in this setting and its differences with the classical case. We show that the compositions of even numbers of vector reflections are not enough to fully describe SO${}_0$ since they only show an orthogonal structure and do not include the symplectic part of SO${}_0$. Next we propose an alternative, by defining the spin group through the exponential of extended superbivectors and show that they indeed cover the whole set of superrotations. In particular, we explicitly describe a subset $\Xi$ which is a double covering of SO${}_0$ and contains in particular every fractional Fourier transform.

\section{Preliminaries}\label{Prem}
In this section we provide some preliminaries on the Grassmann algebra of coefficients that is going to be used to define linear actions on supervector variables. The definitions and properties contained in this section are well known material and can be found for example in \cite{MR1701597, MR2069561, Berezin:1987:ISA:38130, MR2840967}. We provide this summary to increase the readability of this manuscript.

\subsection{Grassmann algebras and Grassmann envelopes}
Let $\Lam_N$ be the Grassmann algebra of order $N\in \N$ over the field $\K$ ($\K=\R$ or $\C$) with canonical generators $f_1, \ldots, f_N$ which are subject to the multiplication rules $f_jf_k+f_kf_j=0$
implying in particular that $f_j^2=0$.  A basis for ${\Lam}_N$ consists of elements of the form $f_{\emptyset}=1$, $ f_A=f_{j_1}\cdots f_{j_k}$ for $A=\{j_1, \dots, j_k\}$ $(1\leq j_1< \ldots< j_k\leq N)$. Hence an arbitrary element $a\in {\Lam}_N$ has the form $a=\sum_{A\inc \{1, \ldots, N\}}a_Af_A$ with $a_A\in\K$.
We define the space of homogeneous elements of degree $k$ by $\Lam_N^{(k)}=\mbox{span}_\K\{f_A:\;|A|=k\}$,
 where in particular $\Lam_N^{(k)}=\{0\} $ for $k>N$. It then easily follows that $\Lam_N=\bigoplus_{k=0}^N \Lam_N^{(k)}$ and $\Lam_N^{(k)} \Lam_N^{(\ell)} \inc \Lam_N^{(k+\ell)}$.

The projection of $\Lam_N$ on its $k$-homogeneous part is denoted by $[\cdot]_k:\Lam_N\fd \Lam_N^{(k)}$, i.e. $[a]_k=\sum_{|A|=k} a_Af_A$. 
 In particular we denote $[a]_0=a_{\emptyset}=:a_0$.
It is well-known that $\Lam_N$ shows a natural $\Z_2$-grading.  In fact, defining $\Lam_N^{(ev)} = \bigoplus_{k\geq0} \Lam_N^{(2k)}$ and $\Lam_N^{(odd)}=\bigoplus_{k\geq0} \Lam_N^{(2k+1)}$
as the spaces of homogeneous even and odd elements respectively, we obtain the superalgebra structure $\Lam_N=\Lam_N^{(ev)}\oplus \Lam_N^{(odd)}$.
We recall that $\Lam_{N}$ is graded commutative in the sense that 
\[vw=wv, \hspace{.5cm} vw\p=w\p v, \hspace{.5cm} v\p w\p=-w\p v\p, \hspace{1cm} v,w\in \Lam_N^{(ev)}, \;\;\; v\p,w\p \in\Lam_N^{(odd)}. \] 
Every $a\in {\Lam}_N$ may be written as the sum $a=a_0+{\bf a}$ of a number $a_0:=a_\emptyset \in \K$ and a nilpotent element ${\bf a}=\sum_{|A|\leq 1}a_Af_A$ (in particular ${\bf a}^{N+1}=0$). The elements $a_0$, ${\bf a}$ are called the body and the nilpotent part of $a\in{\Lam}_N$, respectively. The subalgebra of all nilpotent elements is denoted by $\Lam_N^+:= \bigoplus_{k=1}^N \Lam_N^{(k)}$.
It is easily seen that the projection $[\cdot]_0:\Lam_N\fd \K$ is an algebra homomorphism, i.e.\ $[ab]_0=a_0b_0$ for $a,b\in \Lam_N$.
In particular the following property holds.
\begin{lem}\label{Sq}
Let $a\in\Lam_N$ such that $a^2\in \K\setminus\{0\}$. Then $a\in\K$.
\end{lem}


The algebra $\Lam_N$ is a $\K$-vector space of dimension $2^N$. As every finite dimensional $\K$-vector space, $\Lam_N$ becomes a Banach space with the introduction of an arbitrary norm, all norms being equivalent. In particular, 
the norm  $\|\cdot\|_{\Lam}$ defined on $\Lam_N$ by $\|a\|_{\Lam}=\sum_{A\inc\{1,\ldots,N\}}|a_A|$ satisfies 
\[
\|ab\|_{\Lam} \leq \|a\|_{\Lam} \|b\|_{\Lam}, \hspace{.5cm} \mbox{ for every \ } a,b\in \Lam_N.
\]
%
The exponential of $a\in \Lam_N$, denoted by $e^a$ or $\exp(a)$, is defined by the power series 
\[e^a=\sum_{j=0}^\infty \frac{a^j}{j!}.\]
This series converges for every $a\in \Lam_N$ and defines a continuous function in $\Lam_N$.
Now consider the graded vector space $\K^{p,q}$ with standard homogeneous basis $e_1, \ldots, e_p, e\p_1, \ldots, e\p_q$, i.e.\ $\K^{p,q}=\K^{p,0}\oplus \K^{0,q}$ where $\{e_1, \ldots, e_p\}$ is a basis for $\K^{p,0}$ and $\{e\p_1, \ldots, e\p_q\}$ is a basis for $\K^{0,q}$. Elements in $\K^{p,0}$ are called {\it even homogeneous elements} while elements in $\K^{0,q}$ are called {\it odd homogeneous elements}. In \cite[p.~91]{Berezin:1987:ISA:38130}, the {\it Grassmann envelope} $\K^{p,q}(\Lam_N)$ was defined as the set of Grassmann supervectors 
\begin{equation}\label{SupVecGrass}
{\bf w}=\underline{w}+\underline{w\p}=\sum_{j=1}^p w_j e_j + \sum_{j=1}^q w\p_je\p_j, \quad \mbox{ where \ } w_j\in \Lam_N^{(ev)}, \;\;\;w\p_j\in \Lam_N^{(odd)}.
\end{equation}
\begin{remark}
The Grassmann envelope of a general graded $\K$-vector space $V= V_{\ba 0} \oplus V_{\ba 1}$ is similarly defined as
\[V(\Lam_N) = \left(V \otimes \Lam_N\right)_0 = \left( \Lam_N^{(ev)} \otimes V_{\ba 0} \right) \oplus \left(\Lam_N^{(odd)} \otimes V_{\ba 1} \right). \]
\end{remark}
The set $\K^{p,q}(\Lam_N)$ is a $\K$-vector space of dimension $2^{N-1}(p+q)$, inheriting the $\Z_2$-grading of $\K^{p,q}$, i.e.\ 
\[\K^{p,q}(\Lam_N)=\K^{p,0}(\Lam_N)\oplus \K^{0,q}(\Lam_N),\]
where $\K^{p,0}(\Lam_N)$ denotes the subspace of vectors of the form (\ref{SupVecGrass}) with $w\p_j=0$, and $\K^{0,q}(\Lam_N)$ denotes the subspace of vectors of the form (\ref{SupVecGrass}) with $w_j=0$. The subspaces $\K^{p,0}(\Lam_N)$ and $\K^{0,q}(\Lam_N)$ are called the Grassmann envelopes of $\K^{p,0}$ and $\K^{0,q}$, respectively.

In $\K^{p,q}(\Lam_N)$, there exists a subspace which is naturally isomorphic to $\K^{p,0}$. It consists of vectors (\ref{SupVecGrass}) of the form ${\bf w}=\sum_{j=1}^m w_je_j$ with $w_j\in \K$. The map $[\cdot]_0: \K^{p,q}(\Lam_N)\fd \K^{p,0}$ defined by $[{\bf w}]_0=\sum_{j=1}^p \left[w_j\right]_0 e_j$ will be useful. 

The standard basis of $\K^{p,q}$ can be represented by the columns $e_j=\left( 0, \ldots, 1, \ldots, 0 \right)^T$ (1 on the $j$-th place from the left) and $e\p_j=\left( 0, \ldots, 1, \ldots, 0 \right)^T$ (1 on the $(p+j)$-th place from the left). In this basis, elements of $\K^{p,q}(\Lam_N)$ take the form ${\bf w}=\left(w_1, \ldots, w_p, w\p_1, \ldots, w\p_q \right)^T$.

\subsection{Supermatrices}
The $\Z_2$-grading of $\K^{p,q}$ yields the $\Z_2$-grading of the space $\mbox{End}\left(\K^{p,q}\right)$ of endomorphisms on $\K^{p,q}$. When seen as a Lie superalgebra, this space is denoted by $\mathfrak{gl}(p|q)(\K)$. The super Lie bracket on $\mathfrak{gl}(p|q)(\K)$ is given by $[X,Y]=XY-(-1)^{|X||Y|} YX$ where $X,Y\in \mathfrak{gl}(p|q)(\K)$ are homogeneous elements, i.e.\ elements in the even or in the odd subalgebra. Here the grading function $|X|$ is defined as $0$ if $X$ is even and $1$ if $X$ is odd. 
%

It is easily seen that the Grassmann envelope of any Lie subsuperalgebra of  $\mathfrak{gl}(p|q)(\K)$ is a classical Lie algebra.
The Grassmann envelope of $\mathfrak{gl}(p|q)(\K)$ is denoted by $\Mat(p|q)(\Lam_N)$. Elements in $\Mat(p|q)(\Lam_N)$ are called supermatrices and are of the form
\begin{equation}\label{SMat}
M=\left(\begin{array}{cc} A & B^{\p} \\ C^{\p} & D\end{array}\right)=\left(\begin{array}{cc} A & 0 \\ 0 & D\end{array}\right)+\left(\begin{array}{cc} 0 & B^{\p} \\ C^{\p} & 0\end{array}\right)
\end{equation}
where $A$ is a $p\times p$ matrix with entries in $\Lam_N^{(ev)}$, $B^{\p}$ is a $p\times q$ matrix with entries in $\Lam_N^{(odd)}$, $C^{\p}$ is a $q\times p$ matrix with entries in $\Lam_N^{(odd)}$, and $D$ is a $q\times q$ matrix with entries in $\Lam_N^{(ev)}$.
Let $\Mat(p|q)(\Lam_N^{(k)})$ be the space of homogeneous supermatrices of degree $k$. These subspaces define a grading in $\Mat(p|q)(\Lam_N)$ by $\Mat(p|q)(\Lam_N)=\bigoplus_{k=0}^N \Mat(p|q)(\Lam_N^{(k)})$ and $\Mat(p|q)(\Lam_N^{(k)}) \; \Mat(p|q)(\Lam_N^{(\ell)}) \inc  \Mat(p|q)(\Lam_N^{(k+\ell)})$.

Every supermatrix $M$ can be written as the sum of a {\it body} matrix $M_0\in  \Mat(p|q)(\Lam_N^{(0)}) $ and a nilpotent element ${\bf M} \hspace{-0.5mm}\in \hspace{-0.5mm}  \Mat(p|q)(\Lam_N^+) \hspace{-0.5mm} := \hspace{-0.5mm} \bigoplus_{k=1}^N \Mat(p|q)(\Lam_N^{(k)})$. 
We also define the algebra homomorphism $[\cdot]_0:\Mat(p|q)(\Lam_N)\fd \Mat(p|q)(\Lam_N^{(0)})$ as the projection:
\[M=\left(\begin{array}{cc} A & B^{\p} \\ C^{\p} & D\end{array}\right) \Fd \left(\begin{array}{cc} A_0 & 0 \\ 0 & D_0\end{array}\right)=M_0=[M]_0\]
where $A_0$  and $D_0$ are the {\it body} projections of $A$ and $D$ on $\K^{p\times p}$ and $\K^{q\times q}$ respectively. We recall that  $\Mat(p|q)(\Lam_N^{(0)})$ is equal to the even subalgebra of $\mathfrak{gl}(p|q)(\K)$. Given a set of supermatrices ${\bf S}$ we define $\left[{\bf S}\right]_0=\{[M]_0:M\in{\bf S}\}$.

Let $\GL(p|q)(\Lam_N)$ be the Lie group of all invertible elements of $\Mat(p|q)(\Lam_N)$. A supermatrix $M \in \Mat(p|q)(\Lam_N)$ of the form (\ref{SMat}) is invertible if and only if the blocks $A$ and $D$ are invertible (or equivalently, $A_0$ and $D_0$ are invertible), see e.g.\ \cite{Berezin:1987:ISA:38130}. The inverse of every $M\in \GL(p|q)(\Lam_N)$ is given by
\begin{equation}\label{Inv}
M^{-1}=\left(\begin{array}{cc} \left(A-B^{\p}D^{-1}C^{\p}\right)^{-1} & -A^{-1}B^{\p} \left(D-C^{\p}A^{-1}B^{\p}\right)^{-1} \\ -D^{-1}C^{\p} \left(A-B^{\p}D^{-1}C^{\p}\right)^{-1} & \left(D-C^{\p}A^{-1}B^{\p}\right)^{-1}\end{array}\right).
\end{equation} 

The {\it supertranspose} of $M\in \Mat(p|q)(\Lam_N)$ is defined by $M^{ST}=\left(\begin{array}{cc} A^T & C^{\p\, T} \\ -B^{\p\, T} & D^T\end{array}\right)$,
where $\cdot^T$ denotes the usual matrix transpose. It easily follows that $\left(M L\right)^{ST}=L^{ST} M^{ST}$.

The {\it supertrace} is defined as the map $\str: \Mat(p|q)(\Lam_N)\fd \Lam_N^{(ev)}$ given by $\mbox{str}(M) = \mbox{tr} (A)- \mbox{tr} (D)$. It follows that $\str (ML)=\str(LM)$ for any pair $M,L\in \Mat(p|q)(\Lam_N)$.

The {\it superdeterminant} or {\it Berezinian} is a function from $\GL(p|q)(\Lam_N)$ to $\Lam_N^{(ev)}$ defined by
\[\sd(M)=\det\left(A-B^{\p}D^{-1}C^{\p}\right)  {\det(D)}^{-1}={\det\left(D-C^{\p}A^{-1}B^{\p}\right)}^{-1} {\det(A)}.\]
In particular, one has that $\sd(ML)=\sd(M)\sd(L)$ and $\sd\left(M^{ST}\right)=\sd(M)$.

In the vector space $\Mat(p|q)(\Lam_N)$ we introduce the norm 
$\|M\|=\sum_{j,k=1}^{p+q}\|m_{j,k}\|_{\Lam},$
 where $m_{j,k}\in \Lam_N$ ($j,k=1,\ldots,p+q$) are the entries of $M\in\Mat(p|q)(\Lam_N)$.
As was the case in $\Lam_N$, also this norm satisfies the inequality $\|ML\|\leq  \|M\|\|L\|$ for every pair $M, L \in \Mat(p|q)(\Lam_N)$, leading to the absolute convergence of the series 
\[\exp(M)=\sum_{j=0}^\infty \frac{M^j}{j!}\]
and hence, the continuity of the exponential map in $\Mat(p|q)(\Lam_N)$. The supertranspose, the supertrace and the superdeterminant are also continuous maps. In particular, one has that (see e.g.\ \cite{Berezin:1987:ISA:38130})
\begin{equation}\label{ExpPro}
\sd\left(e^M\right)=e^{\str(M)}, \;\;\;\;\;\;\;\;\;\;\;\; M\in\Mat(p|q)(\Lam_N). 
\end{equation}
Moreover, $e^{tM}$ ($t\in \R$) is a smooth curve  with $\frac{d}{dt}e^{tM}=M e^{tM}=e^{tM}M$ and $\frac{d}{dt}e^{tM}\big|_{t=0}=M$.

The logarithm for a supermatrix $M\in \Mat(p|q)(\Lam_N)$ is defined by $\ln(M)=\sum_{j=1}^\infty (-1)^{j+1}\frac{(M-I_{p+q})^j}{j}$
wherever  it converges. This series onverges and yields a continuous function near $I_{p+q}$. 
\begin{pro}
In $\Mat(p|q)(\Lam_N)$, let $U$ be a neighbourhood of $I_{p+q}$ on which $\ln$ is defined and let $V$ be a neighbourhood of $0$ such that $\exp(V):=\{e^M:M\in V\}\inc U$. Then $e^{\ln (M)}=M$, $\pt M\in U$, and $\ln (e^L)=L$, $\pt L\in V$.
\end{pro}
The exponential of a nilpotent matrix ${\bf M}\in \Mat(p|q)(\Lam_N^+)$ reduces to a finite sum, yielding the bijective mapping
\[\exp:\Mat(p|q)(\Lam_N^+)\fd I_{p+q}+ \Mat(p|q)(\Lam_N^+)\]
with inverse
\[\ln: I_{p+q}+ \Mat(p|q)(\Lam_N^+) \fd \Mat(p|q)(\Lam_N^+),\]
since also the second expansion only has a finite number of non-zero terms, whence problems of convergence do not arise. We recall that a supermatrix $M$ belongs to $\GL(p|q)(\Lam_N)$ if and only if its body $M_0$ has an inverse. Then 
$M=M_0(I_{p+q}+M_0^{-1}{\bf M})=M_0\exp({\bf L}), 
$
for some unique ${\bf L}\in\Mat(p|q)(\Lam_N^+)$.


\section{The algebra $\mathcal{A}_{m,2n}\otimes \Lam_N$}\label{SsFwGC}
\subsection{The Clifford setting in superspace}
In order to set up the Clifford analysis framework in superspace, take $p=m$, $q=2n$ ($m,n\in\N$) and $\K=\R$. The canonical homogeneous basis $e_1,\ldots, e_m, e\p_1, \ldots, e\p_{2n}$ of $\R^{m,2n}$ can be endowed with an orthogonal and a symplectic structure by the multiplication rules
\begin{equation}\label{Reg Gen}
e_je_k+e_ke_j=-2\del_{j,k}, \hspace{.3cm} e_je\p_k+e\p_ke_j=0, \hspace{.3cm} e\p_je\p_k-e\p_ke\p_j=g_{j,k},
\end{equation}
where the symplectic form $g_{j,k}$ is defined by
\[g_{2j,2k}=g_{2j-1,2k-1}=0, \hspace{.5cm} g_{2j-1,2k}=-g_{2k,2j-1}=\del_{j,k}, \hspace{.5cm} j,k=1,\ldots,n.\]
Following these relations, elements in $\R^{m,2n}$ generate an infinite dimensional algebra denoted by $\mathcal C_{m,2n}$.

{The definition of the Clifford supervector variable follows from a representation of the so-called radial algebra, see e.g.\ \cite{273847, Sommen_2002, Bie2007, MR3656057}. Given a set $S=\{x, y, \ldots\}$ of $\ell>1$ abstract vector variables we define the {\it radial algebra} $R(S)$ as the associative algebra over $\R$ freely generated by $S$ and subject to the defining axiom
\[
\mbox{\bf (A1)} \hspace{.5cm} \left[\{x,y\},z\right]=0 \hspace{.3cm} \mbox{ for any } x,y,z\in S,
\]
where $\{a,b\}=ab+ba$ and $[a,b]=ab-ba$. A radial algebra  representation is an algebra homomorphism $\Psi:R(S)\fd \mathfrak{A}$ from $R(S)$ into an algebra $\mathfrak{A}$. The term representation also refers to the range $\Psi(R(S))\inc \mathfrak{A}$ of that mapping. The easiest and most important example of radial algebra representation is the algebra generated by standard Clifford vector variables. In that way, $R(S)$ describes the main algebraic properties of the Clifford function theory. For a detailed study on radial algebras and its representations we refer the reader to \cite{273847, Sommen_2002, Bie2007, MR3656057}.}


The representation of the radial algebra in superspace is defined by the mapping 
\begin{equation}\label{SupRep}
x\mapsto {\bf x}=\underline{x}+\underline{x}\p=\sum_{j=1}^m x_j e_j+\sum_{j=1}^{2n}x\p_je\p_j, \hspace{.3cm} x\in S,
\end{equation}
between $S$ and the set of independent supervector variables ${\bf S}=\{{\bf x}: x\in S\}$. For each $x\in S$ we consider in (\ref{SupRep}) $m$ bosonic (commuting) variables $x_1,\ldots, x_m$ and $2n$ fermionic (anti-commuting) variables $x\p_1,\ldots, x\p_{2n}$.  The projections $\underline{x}=\sum_{j=1}^m x_j e_j$ and $\underline{x}\p=\sum_{j=1}^{2n}x\p_j e\p_j$ are called the bosonic and fermionic vector variables, respectively.

Let us define  the sets $VAR$ and $VAR^{\,}\p$ of bosonic and fermionic variables 
\begin{equation}\label{VArVar_prime}
VAR=\bigcup_{{\bf x}\in {\bf S}}\{x_1,\ldots, x_m\}, \hspace{1cm} VAR^{\,}\p=\bigcup_{{\bf x}\in {\bf S}}\{x\p_1,\ldots, x\p_{2n}\}
\end{equation}
respectively, where the sets $\{x_1,\ldots, x_m\}$ and $\{x\p_1,\ldots, x\p_{2n}\}$ correspond to the coordinates of the bosonic and the fermionic vector variables associated to each ${{\bf x}\in {\bf S}}$ through the correspondences (\ref{SupRep}). In this way, $VAR$ contains $m\ell$ bosonic variables and $VAR^{\,}\p$ contains $2n\ell$ fermionic variables. They give rise to the algebra of super-polynomials $\mathcal V=\mbox{Alg}_{\R}\{VAR\cup VAR^{\,}\p\}$ which is extended to the algebra of Clifford-valued super-polynomials 
\[\mathcal A_{m,2n}=\mathcal V \otimes \mathcal C_{m,2n},\]
where the elements of $\mathcal{V}$ commute with the elements of  $\mathcal C_{m,2n}$.

The algebra $\mathcal V$ clearly is  $\Z_2$-graded. Indeed,  $\mathcal V={\mathcal V}_{\ba 0}\oplus {\mathcal V}_{\ba 1}$ where ${\mathcal V}_{\ba 0}$ consists of all commuting super-polynomials and ${\mathcal V}_{\ba 1}$ consists of all anti-commuting super-polynomials in $\mathcal V$. It is easily seen that the fundamental axiom of the radial algebra is fulfilled in this representation since for every pair ${\bf x},{\bf y}\in {\bf S}$
\begin{equation}\label{A-Com}
\langle {\bf x},{\bf y}\rangle := -\frac{1}{2}\{{\bf x},{\bf y}\}=\sum_{j=1}^m x_jy_j-\frac{1}{2}\sum_{j=1}^n(x\p_{2j-1}y\p_{2j}-x\p_{2j}y\p_{2j-1})\in  \mathcal{V}_{\ba 0}.
\end{equation}
The above definition $\langle {\bf x},{\bf y}\rangle$ will be used as generalized inner product. As mentioned before, one of the goals of this paper is to study the invariance under linear transformations of this inner product of supervector variables. In order to study linear actions on the algebra $\mathcal{A}_{m,2n}=\mathcal{V}\otimes \mathcal{C}_{m,2n}$ we must consider a suitable set of coefficients. Observe that the field of numbers $\K=\R$ or $\C$ is too limited for that aim since it does not lead to any interaction between even and odd elements. For instance, multiplication by real or complex numbers leaves the  decomposition $\mathcal V= \mathcal V_{\ba 0} \oplus \mathcal V_{\ba 1}$ of the algebra of super-polynomials invariant.

The study of linear actions on $\mathcal{A}_{m,2n}$ requires of a set including both commuting and anti-commuting elements. In this paper we consider the most simple set of such coefficients, i.e.\ the Grassmann algebra $\Lam_N$ generated by odd independent elements $f_1\ldots, f_N$.  This leads to the $\Z_2$-graded algebra of super-polynomials with Grassmann coefficients 
$\mathcal{V} \otimes \Lam_N
,$
generated over $\R$ by the set of $m\ell$ commuting variables $VAR$ and the set of independent $2n\ell+N$ anti-commuting symbols $VAR^{\,}\p \, \cup\{f_1, \ldots, f_N\}$. In general we consider the algebra 
\[\mathcal{A}_{m,2n}\otimes \Lam_N=  \mathcal{V} \otimes \Lam_N \otimes \mathcal{C}_{m,2n},\]
of super-polynomials with coefficients in $\Lam_N \otimes \mathcal{C}_{m,2n}$. Here elements of $ \mathcal{V} \otimes \Lam_N$ commute with elements in $\mathcal{C}_{m,2n}$. 

In the set of coefficients $\Lam_N \otimes \mathcal{C}_{m,2n}$ one has a radial algebra representation by considering supervectors ${\bf w}\in \R^{m,2n}(\Lam_N)$, i.e.\
\[{\bf w}=\underline{w}+\underline{w\p}=\sum_{j=1}^m w_je_j+\sum_{j=1}^{2n} w\p_je\p_j, \hspace{1cm} w_j\in\Lam_N^{(ev)}, \;\;\;w\p_j\in\Lam_N^{(odd)}, \]
where clearly the basis elements $e_1,\ldots, e_m, e\p_1, \ldots, e\p_{2n}$ of $\R^{m,2n}$ have to satisfy to the multiplication rules (\ref{Reg Gen}). Indeed, the anti-commutator of two constant supervectors ${\bf w}, {\bf v}\in \R^{m,2n}(\Lam_N)$ clearly is a central element in $\Lam_N \otimes \mathcal{C}_{m,2n}$, i.e.\
$\{ {\bf w}, {\bf v} \}= 2\sum_{j=1}^m w_jv_j+\sum_{j=1}^n(w\p_{2j-1}v\p_{2j}-w\p_{2j}v\p_{2j-1})\in \Lam_N^{(ev)}.$

The subalgebra generated by the Grassmann envelope $\R^{m,2n}(\Lam_N)$ of constant supervectors is called the {\it radial algebra embedded in } $\Lam_N \otimes \mathcal{C}_{m,2n}$. This algebra is denoted by $\R_{m|2n}(\Lam_N)$. Observe that $\R_{m|2n}(\Lam_N)$ is a finite dimensional vector space since it is generated by  the set
\[\left\{f_Ae_j\;:\;A\inc\{1,\ldots,N\}, |A| \mbox{ even }, j=1,\ldots,m\right\} \cup \left\{f_Ae\p_j\;:\;A\inc\{1,\ldots,N\}, |A| \mbox{ odd }, j=1,\ldots,2n\right\},\]
and there is a finite number of possible products amongst these generators.

Every element in $\Lam_N \otimes \mathcal{C}_{m,2n}$ can be written as a finite sum of terms of the form $ae_{j_1}\cdots e_{j_k} e\p_1^{\,\al_1}\cdots e\p_{2n}^{\,\al_{2n}}$ where $a\in\Lam_N$, $1\leq j_1\leq\ldots\leq j_k\leq m$ and $(\al_1,\ldots,\al_{2n})\in \left(\N\cup\{0\}\right)^{2n}$ is a multi-index. In this algebra we consider the corresponding generalization of the projection $[\cdot]_0$ which now goes from $\Lam_N \otimes \mathcal{C}_{m,2n}$ to $\mathcal C_{m,2n}$ and is defined by $\left[ae_{j_1}\cdots e_{j_k} e\p_1^{\,\al_1}\cdots e\p_{2n}^{\,\al_{2n}}\right]_0=[a]_0e_{j_1}\cdots e_{j_k} e\p_1^{\,\al_1}\cdots e\p_{2n}^{\,\al_{2n}}$.

We now can define linear actions on supervector variables ${\bf x}\in {\bf S}$ by means of supermatrices $M\in \Mat(m|2n)(\Lam_N)$. We recall that the basis elements $e_1,\ldots, e_m, e\p_1, \ldots, e\p_{2n}$ can be written as column vectors. Then, by writing the ${\bf x}=\underline{x}+\underline{x\p}\in {\bf S}$ in its column representation we obtain, 
\begin{equation}\label{SupMatActonSupVecVar}
M{\bf x}=\left(\begin{array}{cc} A & B^{\p} \\ C^{\p} & D\end{array}\right)\left(\begin{array}{c}  \underline{x} \\ \underline{x\p}\end{array}\right)= \left(\begin{array}{c}  A\underline{x}+B^{\p}\underline{x\p} \\ C^{\p}\underline{x}+ D\underline{x\p}\end{array}\right).
\end{equation}
This action produces a new supervector variable $M{\bf x}=(y_1,\dots, y_m,y\p_1,\ldots, y\p_{2n})^{T}$ where the $y_j$ are even elements of $\mathcal{V}\otimes \Lam_N$ while the $y\p_j$ are odd ones. It is clear that $(M{\bf x})^{T}={\bf x}^T M^{ST}.$

\subsection{Superbivectors}\label{SupBi-VecExtSect}
Superbivectors in $\Lam_N \otimes \mathcal{C}_{m,2n}$ play a very important r\^ole when studying the invariance of the inner product (\ref{A-Com}). Following the radial algebra approach, the space of bivectors is generated by the wedge product of supervectors of $\R^{m,2n}(\Lam_N)$, i.e.\
\begin{multline*}
\hspace{-.3cm}{\bf w}\wedge {\bf v} =\frac{1}{2}[{\bf w}, {\bf v}]
= \hspace{-.2cm} \sum_{1\leq j<k\leq m}(w_jv_k-w_kv_j)e_je_k+\hspace{-.2cm}\sum_{1\le j\le m\atop 1\le k\le 2n}(w_jv\p_k-w\p_kv_j)e_je\p_k+\hspace{-.2cm}\sum_{1\leq j\leq k\leq 2n}(w\p_jv\p_k+w\p_kv\p_j) \,e\p_j\odot e\p_k,
\end{multline*}
where $e\p_j\odot e\p_k = \frac{1}{2}\{e\p_j,e\p_k\}$. Hence, the space $\R_{m|2n}^{(2)}(\Lam_N)$ of superbivectors consists of elements of the form
\begin{equation}\label{SBiv}
B=\sum_{1\leq j<k\leq m}b_{j,k}\,e_je_k+\sum_{1\le j\le m\atop 1\le k\le 2n}b\p_{j,k}\,e_je\p_k+\sum_{1\leq j\leq k\leq 2n} B_{j,k} \,e\p_j\odot e\p_k, 
\end{equation}
where $b_{j,k}\in \Lam_N^{(ev)}$,  $b\p_{j,k}\in \Lam_N^{(odd)}$ and $B_{j,k}\in \Lam_N^{(ev)}\cap \Lam_N^+$. Observe that the coefficients $B_{j,k}$ are commuting but nilpotent since they are generated by elements of the form $w\p_jv\p_k+w\p_kv\p_j$ that belong to $\Lam_N^+$. This constitutes an important limitation for the space of superbivectors because it means that $\R_{m|2n}^{(2)}(\Lam_N)$ does not allow for any other structure than the orthogonal one. In fact, the real projection $[B]_0$ of every superbivector $B$ is just the classical Clifford bivector: 
\[[B]_0=\sum_{1\leq j<k\leq m}\left[b_{j,k}\right]_0\,e_je_k \in \R_{0,m}^{(2)}.\]
Hence it is necessary to introduce an extension $\R_{m|2n}^{(2)E}(\Lam_N)$ of $\R_{m|2n}^{(2)}(\Lam_N)$  containing elements $B$ of the form (\ref{SBiv}) but with $B_{j,k}\in \Lam_N^{(ev)}$. This extension enables us to consider two different structures in the same element $B$: the orthogonal and the symplectic one. In fact, in this case we have
\[[B]_0=\sum_{1\leq j<k\leq m}\left[b_{j,k}\right]_0\,e_je_k + \sum_{1\leq j\leq k\leq 2n} \left[B_{j,k}\right]_0 \,e\p_j\odot e\p_k.\]
\begin{remark}\label{DimBiv}
Observe that $\R_{m|2n}^{(2)}(\Lam_N)$ and  $\R_{m|2n}^{(2)E}(\Lam_N)$ are finite dimensional real vector subspaces of $\Lam_N \otimes \mathcal{C}_{m,2n}$ with
\begin{align*}
\dim \R_{m|2n}^{(2)}(\Lam_N) &=  2^{N-1}\frac{m(m-1)}{2}+2^{N-1}\,2mn+\left(2^{N-1}-1\right)\,n(2n+1),\\
\dim \R_{m|2n}^{(2)E}(\Lam_N)&= 2^{N-1}\frac{m(m-1)}{2}+2^{N-1}\,2mn+2^{N-1}\,n(2n+1).
\end{align*}
\end{remark} 
The extension $\R_{m|2n}^{(2)E}(\Lam_N)$ of the superbivector space clearly lies outside the radial algebra $ \R_{m|2n}(\Lam_N)$, and its elements generate an infinite dimensional algebra. 
Elements in $\R_{m|2n}^{(2)E}(\Lam_N)$ are called {\it extended superbivectors}. Both superbivectors and extended superbivectors preserve several properties of classical Clifford bivectors.
\begin{pro}
The space $\R_{m|2n}^{(2)E}(\Lam_N)$ is a Lie algebra. In addition, $\R_{m|2n}^{(2)}(\Lam_N)$ is a Lie subalgebra of $\R_{m|2n}^{(2)E}(\Lam_N)$.
\end{pro}
\pf
We only need to check that the Lie bracket defined by the commutator in the associative algebra $\mathcal{A}_{m,2n}\otimes \Lam_N$ is an internal binary operation in $\R_{m|2n}^{(2)E}(\Lam_N)$  and $\R_{m|2n}^{(2)}(\Lam_N)$. Direct computation shows that for $a,b\in \Lam_N^{(ev)}$ and $a\p,b\p\in \Lam_N^{(odd)}$ we get:
\begin{align*}
[ae_je_k, \,be_re_s] &= ab\,(2\del_{j,s}\,e_re_k-2\del_{s,k}\,e_re_j+2\del_{r,j}\,e_ke_s-2\del_{r,k}\,e_je_s),\\
[ae_je_k, \,b\p e_re\p_s] &= ab\p\,(2\del_{r,j}\,e_ke\p_s-2\del_{r,k}\,e_je\p_s),\\
[ae_je_k, \,be\p_r\odot e\p_s] &= 0,\\
[a\p e_je\p_k, \,b\p e_re\p_s] &= a\p b\p\;(2\del_{r,j}\,e\p_k \odot e\p_s+(1-\del_{j,r})g_{s,k}\,e_je_r),\\
[a\p e_je\p_k, \,be\p_r\odot e\p_s]  &= a\p b\,(g_{k,s}\,e_je\p_r + g_{k,r}\,e_je\p_s),\\
[ae\p_j\odot e\p_k, \,be\p_r\odot e\p_s]  &= ab\,(g_{j,s}\,e\p_r\odot e\p_k +g_{k,s}\,e\p_r\odot e\p_j + g_{j,r}\,e\p_k\odot e\p_s + g_{k,r}\,e\p_j\odot e\p_s).
\end{align*}
$\hfill\square$

It is well known from the radial algebra framework that the commutator of a bivector with a vector always yields a linear combination of vectors with coefficients in the scalar subalgebra. Indeed, for the abstract vector variables $x,y,z\in S$ we obtain
\[[x\wedge y, z]=\frac{1}{2}\left[[x,y],z\right]=\frac{1}{2}[2xy-\{x,y\}, z]=[xy,z]=\{y,z\}x-\{x,z\}y.\]
This property can be easily generalized to $\R_{m|2n}^{(2)E}(\Lam_N)$ by straightforward computation. In particular, the following results hold. 
\begin{pro}\label{ComBivVec}
Let ${\bf x}\in {\bf S}$ be a supervector variable, let $\{b_1, \ldots, b_{2^{N-1}}\}$ be a basis for $\Lam_N^{(ev)}$ and let $\{b\p_1, \ldots, b\p_{2^{N-1}}\}$ be a basis for $\Lam_N^{(odd)}$. Then,
\begin{align*}
[b_r\,e_je_k,\,{\bf x}]  &= 2b_r\,(x_je_k-x_ke_j),                                                   &                  [b_r \,e\p_{2j}\odot e\p_{2k},\,{\bf x}]  &= -b_r\,(x\p_{2j-1}e\p_{2k}+x\p_{2k-1}e\p_{2j}), \\ 
[b\p_r\,e_je\p_{2k-1},\,{\bf x}]  &= b\p_r\,(2x_je\p_{2k-1}+x\p_{2k}e_j),                 &                  [b_r \,e\p_{2j-1}\odot e\p_{2k-1},\,{\bf x}]  &= b_r\,(x\p_{2j}e\p_{2k-1}+x\p_{2k}e\p_{2j-1}),\\
[b\p_r\,e_je\p_{2k},\,{\bf x}]  &= b\p_r\,(2x_je\p_{2k}-x\p_{2k-1}e_j),                    &                  [b_r \,e\p_{2j-1}\odot e\p_{2k},\,{\bf x}]  &= b_r\,(x\p_{2j}e\p_{2k}-x\p_{2k-1}e\p_{2j-1}).
\end{align*}
\end{pro}
\noindent Clearly, the above computations remain valid when replacing ${\bf x}$ by a fixed supervector ${\bf w}\in \R^{m,2n}(\Lam_N)$.

\subsection{Tensor algebra and exponential map}
Since  $\Lam_N \otimes \mathcal{C}_{m,2n}$ is infinite dimensional, the definition of the exponential map by means of a power series 
is not as straightforward as it was for the algebras $\Lam_N$ or $\Mat(p|q)(\Lam_N)$. A correct definition of the exponential map in $\Lam_N \otimes \mathcal{C}_{m,2n}$ requires the introduction of the tensor algebra. More details about the general theory of tensor algebras can be found in several basic references, see e.g.\ \cite{MR1187759, MR1920389, DIENES01011930}.

{Let $V$ be the $\R$-vector space with basis $B_V=\{f_1,\ldots,f_N,e_1,\ldots,e_m,e\p_1,\ldots,e\p_{2n}\}$ and let $T(V)$ be its tensor algebra. This is $T(V)=\bigoplus_{j=0}^\infty T^j(V)$ where $T^j(V)=\mbox{span}_\R\{v_1\otimes\cdots\otimes v_j : v_\ell\in B_V\}$ is the $j$-fold tensor product of $V$ with itself. }
Then $\Lam_N \otimes \mathcal{C}_{m,2n}$ can be seen as a subalgebra of $T(V)/I$ where $I\inc T(V)$ is the two-sided ideal generated by the elements:
\begin{align*}
f_j\otimes f_k+f_k\otimes f_j,         & &  &f_j\otimes e_k-e_k\otimes f_j,  &&  & f_j\otimes e\p_k-e\p_k\otimes f_j, \\
e_j\otimes e_k+e_k\otimes e_j+2\del_{j,k}, & & &e_j\otimes e\p_k+e\p_k\otimes e_j, & & &e\p_j\otimes e\p_k-e\p_k\otimes e\p_j-g_{j,k}.
\end{align*}  
Indeed, $T(V)/I$ is isomorphic to the extension of $\Lam_N \otimes \mathcal{C}_{m,2n}$ which also contains {\bf infinite} sums of arbitrary terms of the form $ae_{j_1}\cdots e_{j_k} e\p_1^{\,\al_1}\cdots e\p_{2n}^{\,\al_{2n}}$ where $a\in\Lam_N$, $1\leq j_1\leq\ldots\leq j_k\leq m$ and $(\al_1,\ldots,\al_{2n})\in \left(\N\cup\{0\}\right)^{2n}$ is a multi-index.

The exponential map $\exp(x)=e^x=\sum_{j=0}^\infty \frac{x^j}{j!}$ is known to be well defined in the tensor algebra $T(V)$, see e.g.\ \cite{DIENES01011930}, whence it also is well defined in $T(V)/I$. {This approach to the exponential map will be enough for the purposes of this paper. For a more detailed treatment on non flat supermanifolds, we refer the reader to \cite{MR3084559}.}
The following mapping properties hold
\[\exp : \Lam_N \otimes \mathcal{C}_{m,2n} \fd T(V)/I, \hspace{.5cm} \exp :  \R_{m|2n}(\Lam_N)\fd \R_{m|2n}(\Lam_N).  \]
The first statement directly follows from the definition of $T(V)/I$, while the second one can be obtained following the standard procedure established for $\Lam_N$ and $\Mat(p|q)(\Lam_N)$, since the radial algebra $\R_{m|2n}(\Lam_N)\inc  \Lam_N \otimes \mathcal{C}_{m,2n} $ is finite dimensional. 

\section{The orthosymplectic structure in $\R^{m,2n}(\Lam_N)$ }\label{InvIPS}
\subsection{Invariance of the inner product}
The inner product (\ref{A-Com}) can be easily written as 
\[\langle {\bf x}, {\bf y}\rangle={\bf x}^T{\bf Q}{\bf y}\]
in terms of the supermatrix ${\bold Q}={\left(\begin{array}{cc} I_m & 0 \\ 0 & -\frac{1}{2}\Om_{2n}\end{array}\right)}$, where $\Om_{2n}=\mbox{diag}\left(\begin{array}{rc} 0 & 1 \\ -1 & 0\end{array} \right)$. In order to find all supermatrices $M\in \Mat(m|2n)(\Lam_N)$ leaving the inner product $\langle\cdot,\cdot\rangle$ invariant, we observe that
\[ \langle M{\bf x},M{\bf y}\rangle=\langle {\bf x},{\bf y}\rangle \iff \left(M{\bf x}\right)^T{\bf Q}M{\bf y}={\bf x}^T{\bf Q}{\bf y}  \iff {\bf x}^T\left(M^{ST}{\bf Q}M-{\bf Q}\right){\bf y}=0,\]
whence the desired set is given by
\[
\mbox{\normalfont O}_0=\mbox{\normalfont O}_0(m|2n)(\Lam_N)=\big\{M\in \Mat(m|2n)(\Lam_N):  M^{ST}{\bold Q}M-{\bold Q}=0\big\},
\]
\begin{remark}\label{InvGrasInProdRem}
It is clear that elements in the above set of supermatrices also leave the same bilinear form in $\R^{m,2n}(\Lam_N)$ invariant, i.e.\ $-\frac{1}{2}\{ M{\bf w}, M{\bf v} \}=-\frac{1}{2}\{ {\bf w}, {\bf v} \}$, for $M\in \mbox{\normalfont O}_0$ and ${\bf w}, {\bf v} \in \R^{m,2n}(\Lam_N)$.
In general, every property that holds for supermatrix actions on supervector variables ${\bf x}\in {\bf S}$ also holds for the same actions on fixed supervectors ${\bf w}\in\R^{m,2n}(\Lam_N)$.
\end{remark}
\begin{teo}\label{OG}
The following statements hold:
\begin{itemize}
\item[(i)]  $\mbox{\normalfont O}_0(m|2n)(\Lam_N)\inc \GL(m|2n)(\Lam_N)$.
\item[(ii)] $\mbox{\normalfont O}_0(m|2n)(\Lam_N)$ is a group under the usual matrix multiplication.
\item[(iii)] $\mbox{\normalfont O}_0(m|2n)(\Lam_N)$ is a closed subgroup of $\GL(m|2n)(\Lam_N)$.
\end{itemize}
Summarizing, $\mbox{\normalfont O}_0(m|2n)(\Lam_N)$ is a Lie group.
\end{teo}
\pf
\begin{itemize}
\item[$(i)$] For every $M\in \mbox{\normalfont O}_0(m|2n)(\Lam_N)$ we have 
\[\left[M^{ST}\right]_0{\bf Q}[M]_0-{\bf Q}=\left[M\right]_0^{T}{\bf Q}[M]_0-{\bf Q}=0, \hspace{.3cm} \mbox{ where \ } \left[M\right]_0=\left(\begin{array}{cc} A_0 & 0 \\ 0 & D_0\end{array}\right).
\]
 This can be rewritten in terms of the real blocks $A_0$ and $D_0$ as $A_0^T A_0=I_m$ and $D_0^T \Om_{2n} D_0=\Om_{2n}$,
implying that $A_0$ and $D_0$ are invertible matrices. Thus $M$ is invertible.
  
 \item[$(ii)$] It suffices to prove that matrix inversion and matrix multiplication are internal operations in  $\mbox{\normalfont O}_0(m|2n)(\Lam_N)$. Both properties follow by straightforward computation. 
 
 \item[$(iii)$] Let $\{M_j\}_{j\in\N}\inc \mbox{\normalfont O}_0(m|2n)(\Lam_N)$ be a sequence that converges to a supermatrix $M\in \Mat(m|2n)(\Lam_N)$. Since algebraic operations are continuous in the space $\Mat(m|2n)(\Lam_N)$ we have
 \[M^{ST}{\bf Q}M-{\bf Q}=\lim_{j\fd\infty} M_j^{ST}{\bf Q}M_j-{\bf Q}=0 ,\]
 which proves the statement.
$\hfill\square$
\end{itemize}

\begin{pro}\label{CharcO_0}
The following statements hold:
\begin{itemize}
\item[(i)] A supermatrix 
$M\in \Mat(m|2n)(\Lam_N)$ of the form (\ref{SMat}) belongs to $\mbox{\normalfont O}_0$ if and only if
\begin{equation}\label{ConCmp}
\begin{cases}
A^TA-\frac{1}{2}C^{\p\, T}\Om_{2n}C^{\p}=I_m,\\
A^TB^{\p}-\frac{1}{2}C^{\p\, T}\Om_{2n}D=0,\\
B^{\p\, T}B^{\p}+\frac{1}{2}D^{T}\Om_{2n}D=\frac{1}{2}\Om_{2n}.
\end{cases}
\end{equation}
\item[(ii)] $\sd(M)=\pm1$ for every $M\in \mbox{\normalfont O}_0$.
\item[(iii)] $\left[\mbox{\normalfont O}_0\right]_0=\mbox{\normalfont O}(m)\times \mbox{\normalfont Sp}_\Om(2n)$.
 \end{itemize}
\end{pro}
\begin{remark}
As usual, {\normalfont O}$(m)$ is the classical orthogonal group in dimension $m$ and {\normalfont Sp}$_\Om(2n)$ is the symplectic group associated to the antisymmetric matrix $\Om_{2n}$, i.e.\  $\mbox{\normalfont Sp}_\Om(2n)=\{D_0\in \R^{2n \times 2n}: D^{T}_0\Om_{2n}D_0=\Om_{2n}\}$.
\end{remark}
\pf
\begin{itemize}
\item[$(i)$] The relation $M^{ST}{\bf Q}M={\bf Q}$ can be written in terms of $A,B^{\p},C^{\p}, D$ as:
\[\left(\begin{array}{cc} A^TA-\frac{1}{2}C^{\p\, T}\Om_{2n}C^{\p} & A^TB^{\p}-\frac{1}{2}C^{\p\, T}\Om_{2n}D \\ -B^{\p\, T}A-\frac{1}{2}D^{T}\Om_{2n}C^{\p} & -B^{\p\, T}B^{\p}-\frac{1}{2}D^{T}\Om_{2n}D\end{array}\right)=\left(\begin{array}{cc} I_m & 0 \\ 0 & -\frac{1}{2}\Om_{2n}\end{array}\right).\]
\item[$(ii)$] The relation $M^{ST}{\bf Q}M={\bf Q}$ implies that $\sd(M)^2\sd({\bf Q})=\sd({\bf Q})$, whence $\sd(M)^2=1$. The statement then follows from Lemma \ref{Sq}.
\item[$(iii)$] See the proof of Theorem \ref{OG} $(i)$.$\hfill\square$
\end{itemize}

\subsection{Group of superrotations SO${}_0$.}\label{SupRot-Grp}
As in the classical way, we now can introduce the set of {\it superrotations} by 
\[\mbox{\normalfont SO}_0=\mbox{\normalfont SO}_0(m|2n)(\Lam_N)=\{M\in \mbox{\normalfont O}_0: \sd(M)=1\}.\]
This is easily seen to be a Lie subgroup of $\mbox{\normalfont O}_0$ with real projection equal to $\mbox{\normalfont SO}(m)\times \mbox{\normalfont Sp}_\Om(2n)$, where $\mbox{\normalfont SO}(m)\inc \mbox{\normalfont O}(m)$ is the special orthogonal group in dimension $m$.
 In fact, the conditions $M^{ST}{\bf Q}M={\bf Q}$ and $\sd(M)=1$ imply that $M_0^{T}{\bf Q}M_0={\bf Q}$ and $\sd(M_0)=1$,
whence 
$$
M_0=\left(\begin{array}{cc} A_0 & 0 \\ 0 & D_0\end{array}\right)
$$ 
with $A_0^TA_0=I_m$, $D_0^T\Om_{2n}D_0=\Om_{2n}$ and $\det(A_0)=\det(D_0)$. But $D_0\in \mbox{\normalfont Sp}_\Om(2n)$ implies $\det(D_0)=1$, yielding $\det(A_0)=1$ and $A_0\in \mbox{\normalfont SO}(m)$.

The following proposition states that, as in the classical case, $\mbox{\normalfont SO}_0$ is connected and in consequence, it is the identity component of $\mbox{\normalfont O}_0$.
\begin{pro}\label{Con}
\textup{SO}${}_0$ is a connected Lie group.
\end{pro}
\pf
Since the real projection $\mbox{\normalfont SO}(m)\times \textup{Sp}_{\Om}(2n)$ of SO${}_0$ is a connected group, it suffices to prove that for every $M\in \textup{SO}_0$ there exist a continuous path inside SO${}_0$ connecting $M$ with  its real projection $M_0$. To that end, let us write $M=\sum_{j=0}^N [M]_j$,  
where $[M]_j$ is the projection of $M$ on $\Mat(m|2n)(\Lam_N^{(j)})$ for each $j=0,1,\ldots, N$. Then, observe that
\begin{align*}
M^{ST}{\bf Q}M={\bf Q} &\iff \left(\sum_{j=0}^N \left[M^{ST}\right]_j\right){\bf Q}\left(\sum_{j=0}^N [M]_j\right)={\bf Q} \ \iff \ \sum_{k=0}^N \left(\sum_{j=0}^k \left[M^{ST}\right]_j {\bf Q} [M]_{k-j}\right)={\bf Q}\\
&\iff M_0^{T}{\bf Q}M_0={\bf Q}, \hspace{.2cm}\mbox{ and } \hspace{.2cm} \sum_{j=0}^k \left[M^{ST}\right]_j {\bf Q} [M]_{k-j}=0, \hspace{.3cm} k=1,\ldots, N.
\end{align*}
Let us now take the path $M(t)=\sum_{j=0}^N t^j[M]_j$. For $t\in[0,1]$ this is a continuous path with $M(0)=M_0$ and $M(1)=M$. In addition, $M(t)^{T}_0{\bf Q}M(t)_0=M_0^T{\bf Q}M_0={\bf Q}$ 
and for every $k=1,\ldots, N$ we have,
\[\sum_{j=0}^k \left[M(t)^{ST}\right]_j {\bf Q} [M(t)]_{k-j}=t^k\sum_{j=0}^k \left[M^{ST}\right]_j {\bf Q} [M]_{k-j}=0.\]
Hence, $M(t)^{ST}{\bf Q}M(t)={\bf Q}$, $t\in[0,1]$. Finally, observe that $\sd(M(t))=1$ for every $t\in[0,1]$, since $\sd(M(t)_0)=\sd(M_0)=1$.
$\hfill\square$

We will now investigate the corresponding Lie algebras of 
$\mbox{\normalfont O}_0$ and $\mbox{\normalfont SO}_0$.
\begin{teo}\label{LieAlg}
\noindent
\begin{itemize}
\item[(i)] The Lie algebra $\mathfrak{so}_0=\mathfrak{so}_0(m|2n)(\Lam_N)$ of  $\mbox{\normalfont O}_0$ coincides with the Lie algebra of $\mbox{\normalfont SO}_0$ and is given by the space of all {\it "super anti-symmetric" supermatrices}
\[\mathfrak{so}_0=\{X\in \Mat(m|2n)(\Lam_N):X^{ST}{\bf Q}+{\bf Q}X=0\}.\]
\item[(ii)] A supermatrix $X\in \Mat(m|2n)(\Lam_N)$ of the form (\ref{SMat})  belongs to $\mathfrak{so}_0$ if and only if
\begin{equation}\label{so0}
\begin{cases}
A^T+A=0,\\
B^{\p}-\frac{1}{2}C^{\p\, T}\Om_{2n}=0,\\
D^{T}\Om_{2n}+\Om_{2n}D=0.
\end{cases}
\end{equation}
\item[(iii)] $\left[\mathfrak{so}_0\right]_0=\mathfrak{so}(m)\oplus\mathfrak{sp}_\Om(2n)$.
\end{itemize}
\end{teo}
\begin{remark}
 As usual, $\mathfrak{so}(m)=\{A_0\in\R^{m\times m}: A_0^T+A_0=0\}$ is the special orthogonal Lie algebra in dimension $m$ and $\mathfrak{sp}_\Om(2n)=\{D_0\in\R^{2n\times 2n}: D_0^{T}\Om_{2n}+\Om_{2n}D_0=0\}$ is the symplectic Lie algebra defined through the antisymmetric matrix $\Om_{2n}$.
\end{remark}
\pf
\begin{itemize}
\item[$(i)$] If $X\in \Mat(m|2n)(\Lam_N)$ is in the Lie algebra of $\mbox{\normalfont O}_0$ then $e^{tX}\in \mbox{\normalfont O}_0$ for every $t\in\R$, i.e.\ $e^{tX^{ST}}{\bf Q}e^{tX}-{\bf Q}=0$. Differentiating at $t=0$ we obtain $X^{ST}{\bf Q}+{\bf Q}X = 0$. On the other hand, if $X\in \Mat(m|2n)(\Lam_N)$ satisfies $X^{ST}{\bf Q}+{\bf Q}X=0$, then $X^{ST}=-{\bf Q}X{\bf Q}^{-1}$. Computing the exponential of $tX^{ST}$ we obtain
\[e^{tX^{ST}}=\sum_{j=0}^\infty \frac{\left({\bf Q}(-tX){\bf Q}^{-1}\right)^j}{j!}
={\bf Q}e^{-tX}{\bf Q}^{-1},\]
which implies that $e^{tX^{ST}}{\bf Q}e^{tX}-{\bf Q}=0$, 
i.e.\ $e^{tX}\in \mbox{\normalfont O}_0$. Then $\mathfrak{so}_0$ is the Lie algebra of $\mbox{\normalfont O}_0$.
 
From (\ref{ExpPro}) it easily follows that the Lie algebra of $\mbox{\normalfont SO}_0$ is given by
\[\{X\in \Mat(m|2n)(\Lam_N):X^{ST}{\bf Q}+{\bf Q}X=0,\; \str(X)=0\}.\]
But $X^{ST}{\bf Q}+{\bf Q}X=0$ implies $\str(X)=0$. In fact, the condition $X^{ST}=-{\bf Q}X{\bf Q}^{-1}$ implies 
\[\str(X^{ST})=-\str({\bf Q}X{\bf Q}^{-1})=-\str(X),\]
yielding $\str(X)=\str(X^{ST})=-\str(X)$ and $\str(X)=0$. Hence, the Lie algebra of $\mbox{\normalfont SO}_0$ is $\mathfrak{so}_0$.

\item[$(ii)$] Observe that the relation $X^{ST}{\bf Q}+{\bf Q}X=0$ can be written in terms of $A, B^{\p}, C^{\p}, D$ as follows:
\[\left(\begin{array}{cc} A^T+A & -\frac{1}{2}C^{\p \,T}\Om_{2n}+B^{\p}\\ -B^{\p\, T} -\frac{1}{2}\Om_{2n}C^{\p} & -\frac{1}{2}D^{T}\Om_{2n}-\frac{1}{2}\Om_{2n}D \end{array}\right)=0.\]

\item[$(iii)$] 
Let $X=\left(\begin{array}{cc} A & B^{\p} \\ C^{\p} & D\end{array}\right)\in \mathfrak{so}_0$, then $X_0=[X]_0=\left(\begin{array}{cc} A_0 & 0 \\ 0 & D_0\end{array}\right)$ 
satisfies $X_0^{ST}{\bf Q}+{\bf Q}X_0=0$.
Using $(ii)$ we obtain $A_0^T+A_0=0$ and $D_0^{T}\Om_{2n}+\Om_{2n}D_0=0$ which implies that $A_0\in\mathfrak{so}(m)$ and $D_0\in\mathfrak{sp}_\Om(2n)$. 
$\hfill\square$
\end{itemize}
\begin{remark}\label{2ndRem_osp}
{The group $\textup{SO}_0$ is isomorphic to the group obtained by the action of the functor of points of $\textup{OSp}(m|2n)$ on $\Lam_N$. We will explicitly describe this isomorphism in terms of the corresponding Lie algebras. Indeed,} the Lie algebra $\mathfrak{so}_0$ constitutes a Grassmann envelope of the orthosymplectic Lie superalgebra $\mathfrak{osp}(m|2n)$. Here we define $\mathfrak{osp}(m|2n)$, in accordance with \cite{MR3375856}, as the subsuperalgebra of $\mathfrak{gl}(m|2n)(\R)$ given by, 
\begin{align*}
\mathfrak{osp}(m|2n) &:=\{X\in \mathfrak{gl}(m|2n)(\R): X^{ST}{\bf G}+{\bf G} X=0\}, & \mbox{with } \;\; {\bf G}= \left(\begin{array}{cc} I_m & 0 \\ 0 & J_{2n}\end{array}\right)  \;\; J_{2n}= \left(\begin{array}{cc} 0 & I_n \\ -I_n & 0\end{array}\right).
\end{align*}
It suffices to note that $\mathfrak{so}_0$ is the Grassmann envelope of
\[\mathfrak{so}_0(m|2n)(\R) :=\{X\in \mathfrak{gl}(m|2n)(\R): X^{ST}{\bf Q}+{\bf Q} X=0\},\]
which is isomorphic to $\mathfrak{osp}(m|2n)$. In order to explicitly find this isomorphism we first need the matrix
\[
R= \left(\begin{array}{cccc|cccc} 
 1 & 0 & \cdots & 0 & 0 & 0 & \cdots & 0 \\
 0 & 0 & \cdots & 0 & 1 & 0 & \cdots & 0 \\
 0 & 1 & \cdots & 0 & 0 & 0 & \cdots & 0 \\
 0 & 0 & \cdots & 0 & 0 & 1 & \cdots & 0 \\
 \vdots & \vdots & \ddots & \vdots & \vdots & \vdots & \ddots & \vdots \\
 0 & 0 & \cdots & 0 & 0 & 0 & \cdots & 0 \\
 0 & 0 & \cdots & 1 & 0 & 0 & \cdots & 1 
 \end{array}\right) 
 \in \mbox{\normalfont{O}}(2n),
\]
which satisfies $R^T J_{2n} R = \Om_{2n}$. Then the mapping $\Phi: \mathfrak{so}_0(m|2n)(\R) \fd \mathfrak{osp}(m|2n)$, given by
\begin{align*}
\Phi (X) &= {\bf R}^{-1} X {\bf R}, & \mbox{ with } && {\bf R}&= \left(\begin{array}{cc} I_m & 0 \\ 0 & i\sqrt{2} R^T\end{array}\right),
\end{align*}
is easily seen to be a Lie superalgebra isomorphism. Indeed, the matrix ${\bf R}$ is such that ${\bf R}^{ST} {\bf Q} {\bf R} = {\bf G}$. As a consequence, for every $X\in \mathfrak{so}_0(m|2n)(\R)$ one has that
\begin{align*}
\Phi (X)^{ST} {\bf G} + {\bf G} \Phi (X) &= {\bf R}^{ST} X^{ST} \left({\bf R}^{-1}\right)^{ST}  {\bf G} + {\bf G}{\bf R}^{-1} X {\bf R}= {\bf R}^{ST} \left( X^{ST}{\bf Q}+{\bf Q} X \right) {\bf R} =0.
\end{align*}
The use of Grassmann envelopes allows to study particular aspects of the theory of Lie superalgebras in terms of classical Lie algebras and Lie groups. The $ \mathfrak{osp}(m|2n)$-invariance of the super Dirac operator $\pa_{\bf X}$ used in \cite{MR3375856}  has been obtained in \cite{DS_Guz_Somm4} in terms of the invariance of $\pa_{\bf x}$ under the action of the Grassmann envelope $\mathfrak{so}_0$ (or equivalently, under the action of the group $\mbox{\normalfont{SO}}_0$).
\end{remark}
The connectedness of $\mbox{\normalfont SO}_0$ allows to write any of its elements as a finite product of exponentials of supermatrices in $\mathfrak{so}_0$, see \cite[p.~71]{MR3331229}. In the classical case, a single exponential suffices for such a description since $\mbox{\normalfont SO}(m)$ is compact and in consequence $\exp:\mathfrak{so}(m)\fd \mbox{\normalfont SO}(m)$ is surjective, see Corollary 11.10 \cite[p.~314]{MR3331229}. This property, however, does not hold in the group of superrotations $\mbox{\normalfont SO}_0$, since the exponential map from $\mathfrak{sp}_\Om(2n)$ to the non-compact Lie group $\mbox{\normalfont Sp}_\Om(2n)\iso \{I_m\}\times \mbox{\normalfont Sp}_\Om(2n) \inc \mbox{\normalfont SO}_0$ is not surjective, whence not every element in $\mbox{\normalfont SO}_0$ can be written as a single exponential of a supermatrix in $\mathfrak{so}_0$. Nevertheless, it is possible to find a decomposition for elements of $\mbox{\normalfont SO}_0$ in terms of a fixed number of exponentials of $\mathfrak{so}_0$ elements.

Every supermatrix $M\in \mbox{\normalfont SO}_0$ has a unique decomposition $M=M_0+{\bf M}=M_0(I_{m+2n}+{\bf L})$ where $M_0$ is its real projection, ${\bf M}\in \Mat(m|2n)(\Lam_N^+)$ its nilpotent projection and ${\bf L}= M_0^{-1}{\bf M}$. We will now separately study the decompositions for $M_0\in \mbox{\normalfont SO}(m)\times \mbox{\normalfont Sp}_\Om(2n)$ and $I_{m+2n}+{\bf L} \in \mbox{\normalfont SO}_0$.

First consider $M_0 \in \mbox{\normalfont SO}(m)\times \mbox{\normalfont Sp}_\Om(2n)$. We already mentioned that $\exp: \mathfrak{so}(m)\fd \mbox{\normalfont SO}(m)$ is surjective, while $\exp: \mathfrak{sp}_\Om(2n) \fd \mbox{\normalfont Sp}_\Om(2n)$ is not. However, it can be proven that 
\[\mbox{\normalfont Sp}_\Om(2n)=\exp(\mathfrak{sp}_\Om(2n))\cdot \exp(\mathfrak{sp}_\Om(2n)),\] 
invoking the following polar decomposition for real algebraic Lie groups, see Proposition 4.3.3 in \cite{MR3025417}.
\begin{pro}
Let $G\inc \GL(p)$ be an algebraic Lie group such that $G=G^T$ and let $\mathfrak{g}$ be its Lie algebra. Then every $A\in G$ can be uniquely written as $A=Re^X$, $R\in G\cap \mbox{\normalfont O}(p)$, $X\in \mathfrak g \cap \Sym(p)$, 
where $\Sym(p)$ is the subspace of all symmetric matrices in $\R^{p\times p}$. 
\end{pro}
\begin{remark}
A subgroup $G\inc \GL(p)$ is called algebraic if there exists a family $\{p_j\}_{j\in \Upsilon}$ of real polynomials
\[p_j(M)=p_j(m_{11}, m_{12}, \ldots, m_{pp})\in \R[m_{11},  \ldots, m_{pp}]\] 
in the entries of the matrix $M\in \R^{p\times p}$ such that
$G=\{M\in \GL(p): p_j(M)=0, \; \pt j\in \Upsilon\}$.
See \cite[p.~73]{MR3025417} for more details. Obviously, the groups $\mbox{\normalfont O}(m)$, $\mbox{\normalfont SO}(m)$, $\mbox{\normalfont Sp}_\Om(2n)$ are algebraic Lie groups.
\end{remark}
Taking $p=2n$ and $G=\mbox{\normalfont Sp}_\Om(2n)$ in the above proposition we get that every symplectic matrix $D_0$ can be uniquely written as $D_0=R_0e^{Z_0}$ with $R_0\in \mbox{\normalfont Sp}_\Om(2n)\cap \mbox{\normalfont O}(2n)$ and $Z_0\in \mathfrak{sp}_\Om(2n)\cap \Sym(2n)$. But the group $\mbox{\normalfont Sp}_\Om(2n)\cap \mbox{\normalfont O}(2n)$ is isomorphic to the unitary group $\mbox{U}(n)=\left\{ L_0\in \C^{n\times n}: \left(L_0^T\right)^c L_0=I_n\right\}$
which is connected and compact. Then the exponential map from the Lie algebra $ \mathfrak{sp}_\Om(2n)\cap \mathfrak{so}(2n)\iso \mathfrak{u}(n)$ is surjective on $\mbox{\normalfont Sp}_\Om(2n)\cap \mbox{\normalfont O}(2n)$ where
$\mathfrak{u}(n)=\left\{ L_0\in \C^{n\times n}: \left(L_0^T\right)^c + L_0=0\right\}$
is the unitary Lie algebra in dimension $n$. This means that $D_0\in \mbox{\normalfont Sp}_\Om(2n)$ can be written as $D_0=e^{Y_0} e^{Z_0}$ with $Y_0\in  \mathfrak{sp}_\Om(2n)\cap \mathfrak{so}(2n)$ and $Z_0\in \mathfrak{sp}_\Om(2n)\cap \Sym(2n)$. Hence, the supermatrix $M_0\in \mbox{\normalfont SO}(m)\times \mbox{\normalfont Sp}_\Om(2n)$ can be decomposed as
\[M_0=\left(\begin{array}{cc} e^{X_0} & 0 \\ 0 & e^{Y_0}e^{Z_0}\end{array}\right)=\left(\begin{array}{cc} e^{X_0} & 0 \\ 0 & e^{Y_0}\end{array}\right)\left(\begin{array}{cc} I_m & 0 \\ 0 & e^{Z_0}\end{array}\right)=e^Xe^Y,\]
where $X=\left(\begin{array}{cc} {X_0} & 0 \\ 0 & {Y_0}\end{array}\right)\in \mathfrak{so}(m)\times [ \mathfrak{sp}_\Om(2n)\cap \mathfrak{so}(2n)]$ and $Y=\left(\begin{array}{cc} 0 & 0 \\ 0 & {Z_0}\end{array}\right)\in\{0_m\}\times[\mathfrak{sp}_\Om(2n)\cap \Sym(2n)]$.

Now consider the element $I_{m+2n}+{\bf L} \in \mbox{\normalfont SO}_0$. As shown at the end of Section \ref{Prem}, the function $\exp:\Mat(m|2n)(\Lam_N^+)\fd I_{m+2n}+ \Mat(m|2n)(\Lam_N^+)$ is a bijection with the logarithmic function 
as its inverse. Then the supermatrix ${\bf Z}= \ln(I_{m+2n}+{\bf L})$ 
satisfies $e^{\bf Z}=I_{m+2n}+{\bf L}$ 
and is nilpotent. Those properties suffice for proving that ${\bf Z}\in \mathfrak{so}_0$. From now on we will denote the set $\mathfrak{so}_0\cap \Mat(m|2n)(\Lam_N^+)$ of nilpotent elements of $\mathfrak{so}_0$ by  $\mathfrak{so}_0(m|2n)(\Lam_N^+)$.
\begin{pro}
Let ${\bf Z}\in \Mat(m|2n)(\Lam_N^+)$ such that $e^{\bf Z}\in \mbox{\normalfont SO}_0$. Then ${\bf Z}\in \mathfrak{so}_0$.
\end{pro}
\pf
Since $e^{{\bf Z}}\in \mbox{\normalfont SO}_0$, it is clear that $e^{t{\bf Z}}\in \mbox{\normalfont SO}_0$ for every $t\in \Z$. Let us prove that the same property holds for every $t\in \R$.  The expression $e^{t{\bf Z}^{ST}}{\bf Q}e^{t{\bf Z}}-{\bf Q}$ 
can be written as the following polynomial in the real variable $t$.
\begin{align*}
P(t) &= e^{t{\bf Z}^{ST}}{\bf Q}e^{t{\bf Z}}-{\bf Q}=\left[\sum_{j=0}^N \frac{t^j({\bf Z}^{ST})^j}{j!}\right]{\bf Q} \left[\sum_{k=0}^N \frac{t^k {\bf Z}^k}{k!}\right]-{\bf Q}\\
&= \sum_{k=1}^N \sum_{j=0}^k \frac{t^j({\bf Z}^{ST})^j}{j!}\,{\bf Q} \,\frac{t^{k-j}{\bf Z}^{k-j}}{(k-j)!}
=\sum_{k=1}^N \frac{t^k}{k!}\left[\sum_{j=0}^k \binom{k}{j} ({\bf Z}^{ST})^j {\bf Q}  {\bf Z}^{k-j} \right] =  \sum_{k=1}^N \frac{t^k}{k!} P_k({\bf Z}),
\end{align*}
where $P_k({\bf Z})= \sum_{j=0}^k \binom{k}{j} ({\bf Z}^{ST})^j {\bf Q}  {\bf Z}^{k-j}$. If $P(t)$ is not identically zero, i.e.\ not all the $P_k({\bf Z})$ are $0$, we can take $k_0\in\{1, 2, \ldots, N\}$ to be the largest subindex for which $P_{k_0}({\bf Z})\neq 0$. Then,
\[\lim_{t\fd\infty} \frac{1}{t^{k_0}} P(t)=\frac{P_{k_0}({\bf Z})}{k_0!} \neq 0,\]
contradicting that $P(\Z)=\{0\}$. So $P(t)$ identically vanishes, yielding $e^{t{\bf Z}}\in \mbox{\normalfont SO}_0$ for every $t\in \R$.$\hfill\square$\\

In this way, we have proven the following result. 
\begin{teo}\label{SO_0 Decomp}
Every supermatrix in $\mbox{\normalfont SO}_0$ can be written as
\[M=e^X e^Y e^{\bf Z},  \hspace{.3cm} \mbox{ with } \hspace{.3cm}
\begin{cases} X\in \mathfrak{so}(m)\times [ \mathfrak{sp}_\Om(2n)\cap \mathfrak{so}(2n)],\\
Y\in \{0_m\}\times[\mathfrak{sp}_\Om(2n)\cap \Sym(2n)], \\
{\bf Z}\in  \mathfrak{so}_0(m|2n)(\Lam_N^+).
\end{cases}
\]
Moreover, the elements  $Y$ and ${\bf Z}$ are unique.\end{teo}

\subsection{Relation with superbivectors.}
Theorem \ref{LieAlg} allows to compute the dimension of $\mathfrak{so}_0$ as a real vector space.

\begin{cor}
The dimension of the real Lie algebra $\mathfrak{so}_0$ is $2^{N-1}\left(\frac{m(m-1)}{2}+2mn+n(2n+1)\right)$.
\end{cor}
\pf
Since $\mathfrak{so}_0$ is the direct sum of the corresponding subspaces of block components $A, B^{\p}, C^{\p}$ and $D$ respectively, it suffices to compute the dimension of each one of them. According to Theorem \ref{LieAlg} (iii) we have:
\begin{align*}
V_1 &= \left\{ \left(\begin{array}{cc} A & 0 \\ 0& 0\end{array}\right) : A^T=-A,\; A\in \left(\Lam_N^{(ev)}\right)^{m\times m} \right\} \iso \Lam_N^{(ev)} \otimes \mathfrak{so}(m), \\
V_2 &= \left\{ \left(\begin{array}{cc} 0 & \frac{1}{2}C^{\p\, T}\Om_{2n} \\ C^{\p} & 0\end{array}\right) : C^{\p} \in \left(\Lam_N^{(odd)}\right)^{2n\times m}\right\}\iso \Lam_N^{(odd)} \otimes \R^{2n\times m},\\
V_3 &= \left\{ \left(\begin{array}{cc} 0 & 0 \\ 0& D\end{array}\right) : D^{T}\Om_{2n}+\Om_{2n}D=0,\; D\in \left(\Lam_N^{(ev)}\right)^{2n\times 2n} \right\} \iso \Lam_N^{(ev)} \otimes \mathfrak{sp}_\Om(2n).\end{align*}
This leads to $\dim V_1=2^{N-1} \frac{m(m-1)}{2}$, $\dim V_2=2^{N-1} m 2n$ and $\dim V_3=2^{N-1} n(2n+1)$.
$\hfill\square$

Comparing this result with the one in Remark \ref{DimBiv} we obtain that $\dim \R_{m|2n}^{(2)E}(\Lam_N)=\dim \mathfrak{so}_0$.  This means that both vector spaces are isomorphic. This isomorphism also holds on the Lie algebra level. Following the classical Clifford approach, the commutator $[B,{\bf x}]$, with $B\in \R_{m|2n}^{(2)E}(\Lam_N)$ and $ \;{\bf x}\in {\bf S}$,
should be the key for the Lie algebra isomorphism.  Proposition \ref{ComBivVec} shows that 
this commutator defines a linear action on the supervector variable ${\bf x}\in {\bf S}$  that can be represented by a supermatrix in $\Mat(m|2n)(\Lam_N)$, see (\ref{SupMatActonSupVecVar}).
 
 \begin{lem}\label{LemIso}
The map  $\phi:\R_{m|2n}^{(2)E}(\Lam_N)\fd  \Mat(m|2n)(\Lam_N)$ defined by
\begin{equation}\label{PhiIso}
 \phi(B) {\bf x}=[B,{\bf x}] \hspace{2cm}  \;B\in \R_{m|2n}^{(2)E}(\Lam_N), \;{\bf x}\in {\bf S},
\end{equation}
takes values in $\mathfrak{so}_0$. In particular, if we consider $\{b_1,\ldots, b_{2^{N-1}}\}$ and $\{b\p_1,\ldots, b\p_{2^{N-1}}\}$ to be the canonical basis of $\Lam_N^{(ev)}$ and $\Lam_N^{(odd)}$ respectively, we obtain the following basis for $\mathfrak{so}_0$.
{\footnotesize 
\begin{align*}
\phi(b_r e_j e_k) &= 2b_r\left(\begin{array}{cc} E_{k,j}-E_{j,k} & 0 \\ 0 & 0\end{array}\right), & 1\leq r\leq 2^{N-1}, & \;\;1\leq j<k\leq m,\\
\phi(b{\p}_r\, e_j e\p_{2k-1}) &= b\p_r\left(\begin{array}{cc} 0& E_{j,2k} \\ 2E_{2k-1,j} & 0\end{array}\right), & 1\leq r\leq 2^{N-1}, & \;\;1\leq j\leq m, \;\; 1\leq k\leq n, \\
\phi(b{\p}_r\, e_j e\p_{2k}) &= b\p_r\left(\begin{array}{cc} 0& -E_{j,2k-1} \\ 2E_{2k,j} & 0\end{array}\right), & 1\leq r\leq 2^{N-1}, & \;\;1\leq j\leq m, \;\; 1\leq k\leq n, \\
\phi(b{}_r\, e\p_{2j} \odot e\p_{2k}) &= -b_r\left(\begin{array}{cc} 0& 0 \\ 0 & E_{2j,2k-1}+E_{2k,2j-1}\end{array}\right), & 1\leq r\leq 2^{N-1}, & \;\; 1\leq j\leq k\leq n, \\
\phi(b{}_r\, e\p_{2j-1} \odot e\p_{2k-1}) &= b_r\left(\begin{array}{cc} 0& 0 \\ 0 & E_{2j-1,2k}+E_{2k-1,2j}\end{array}\right), & 1\leq r\leq 2^{N-1}, & \;\; 1\leq j\leq k\leq n,\\
\phi(b{}_r\, e\p_{2j-1} \odot e\p_{2k}) &= b_r\left(\begin{array}{cc} 0& 0 \\ 0 & E_{2k,2j}-E_{2j-1,2k-1}\end{array}\right), & 1\leq r\leq 2^{N-1}, & \;\; 1\leq j\leq k\leq n,\\
\phi(b{}_r\, e\p_{2j} \odot e\p_{2k-1}) &= b_r\left(\begin{array}{cc} 0& 0 \\ 0 & E_{2j,2k}-E_{2k-1,2j-1}\end{array}\right), & 1\leq r\leq 2^{N-1}, & \;\; 1\leq j<k\leq n,
\end{align*}
}
where $E_{j,k}$ denotes the matrix in which only the element on the crossing of the $j$-th row and the $k$-th column equals 1 and all the other entries are zero. The order of $E_{j,k}$ should be deduced from the context.
\end{lem}
\pf
The above equalities can be directly obtained from Proposition \ref{ComBivVec}, whence we should only check that all supermatrices obtained  above form a basis for $\mathfrak{so}_0$. The matrices $E_{j,k}$ satisfy the relations
\begin{align*}
E_{j,k}^T&=E_{k,j}, &  E_{j,2k-1}J_{2n}&=E_{j,2k}, & E_{j,2k}J_{2n}&=-E_{j,2k-1}, & J_{2n}E_{2j,k}&=E_{2j-1,k}, & J_{2n}E_{2j-1,k}&=-E_{2j,k}. 
\end{align*}
Then
\begin{itemize}
\item for $\phi(b_r e_j e_k)$ we have $A=2b_r\left(E_{k,j}-E_{j,k}\right)$, $B^{\p}=0$, $C^{\p}=0$ and $D=0$, whence
\[A^T=2b_r\left(E_{j,k}-E_{k,j}\right)=-A;\]

\item for $\phi(b{\p}_r\, e_j e\p_{2k-1})$ we have $A=0$, $B^{\p}=b\p_rE_{j,2k}$, $C^{\p}=2b\p_r E_{2k-1,j}$ and $D=0$, whence
\[\frac{1}{2}C^{\p\,T}\Om_{2n}=b\p_r E_{j,2k-1}\Om_{2n}=b\p_r E_{j,2k}=B^{\p};\]

\item for $\phi(b{\p}_r\, e_j e\p_{2k})$ we have $A=0$, $B^{\p}=-b\p_rE_{j,2k-1}$, $C^{\p}=2b\p_r E_{2k,j}$ and $D=0$, whence
\[\frac{1}{2}C^{\p\,T}\Om_{2n}=b\p_r E_{j,2k}\Om_{2n}=-b\p_r E_{j,2k-1}=B^{\p};\]

\item for $\phi(b{}_r\, e\p_{2j} \odot e\p_{2k})$ we have $A=0$, $B^{\p}=0$, $C^{\p}=0$ and $D=-b_r\left(E_{2j,2k-1}+E_{2k,2j-1}\right)$, whence
\begin{align*}
D^T\Om_{2n}+\Om_{2n}D &= -b_r\left(E_{2k-1,2j}\Om_{2n}+E_{2j-1,2k}\Om_{2n}+\Om_{2n}E_{2j,2k-1}+\Om_{2n}E_{2k,2j-1}\right)\\
&= -b_r\left(-E_{2k-1,2j-1}-E_{2j-1,2k-1}+E_{2j-1,2k-1}+E_{2k-1,2j-1}\right)=0;
\end{align*}

\item for $\phi(b{}_r\, e\p_{2j-1} \odot e\p_{2k-1})$ we have $A=0$, $B^{\p}=0$, $C^{\p}=0$ and $D=b_r\left(E_{2j-1,2k}+E_{2k-1,2j}\right)$, whence
\begin{align*}
D^T\Om_{2n}+\Om_{2n}D &= b_r\left(E_{2k,2j-1}\Om_{2n}+E_{2j,2k-1}\Om_{2n}+\Om_{2n}E_{2j-1,2k}+\Om_{2n}E_{2k-1,2j}\right)\\
&= b_r\left(E_{2k,2j}+E_{2j,2k}-E_{2j,2k}-E_{2k,2j}\right)=0;
\end{align*}

\item for $\phi(b{}_r\, e\p_{2j-1} \odot e\p_{2k})$ we have $A=0$, $B^{\p}=0$, $C^{\p}=0$ and $D=b_r\left(E_{2k,2j}-E_{2j-1,2k-1}\right)$, whence
\begin{align*}
D^T\Om_{2n}+\Om_{2n}D &= b_r\left(E_{2j,2k}\Om_{2n}-E_{2k-1,2j-1}\Om_{2n}+\Om_{2n}E_{2k,2j}-\Om_{2n}E_{2j-1,2k-1}\right)\\
&= b_r\left(-E_{2j,2k-1}-E_{2k-1,2j}+E_{2k-1,2j}+E_{2j,2k-1}\right)=0.
\end{align*}
\end{itemize}
The above computations show that all supermatrices obtained belong to $\mathfrak{so}_0$. Direct verification shows that they form a set of $2^{N-1}\frac{m(m-1)}{2}+2^{N-1}\,2mn+2^{N-1}\,n(2n+1)$ linearly independent elements, i.e.\ a basis of  $\mathfrak{so}_0$.
$\hfill\square$ 

\begin{teo}
The map $\phi:\R_{m|2n}^{(2)E}(\Lam_N)\fd \mathfrak{so}_0$ defined in (\ref{PhiIso}) is a Lie algebra isomorphism.
\end{teo}
\pf 
From Lemma \ref{LemIso} it follows that $\phi$ is a vector space isomorphism. In addition, due to the Jacobi identity in the associative algebra $\mathcal{A}_{m,2n}\otimes \Lam_N$ we have for all $B_1, B_2 \in \R_{m|2n}^{(2)E}(\Lam_N)$ and ${\bf x}\in {\bf S}$ that
\begin{align*}
\left[\phi(B_1), \phi(B_2)\right]{\bf x} &= 
 \left[B_1,\left[B_2,{\bf x}\right]\right]+\left[B_2,\left[{\bf x},B_1\right]\right]= \ \left[\left[B_1,B_2\right],{\bf x}\right] = \ \phi\left([B_1,B_2]\right){\bf x},
\end{align*}
yielding $\left[\phi(B_1), \phi(B_2)\right]= \phi\left([B_1,B_2] \right)$, i.e., $\phi$ is a Lie algebra isomorphism.
$\hfill\square$
\begin{remark}
In virtue of Remark \ref{2ndRem_osp}, the algebra $\R_{m|2n}^{(2)E}(\Lam_N)$ is a Grassmann envelope of $\mathfrak{osp}(m|2n)$.
\end{remark}

\section{The Spin group associated to $\mbox{\normalfont SO}_0$}\label{SecSpGSP}
So far we have seen that the Lie algebra $\mathfrak{so}_0$ of the Lie group of superrotations $\mbox{\normalfont SO}_0$ has a realization in $\Lam_N \otimes \mathcal{C}_{m,2n}$ as the Lie algebra of extended superbivectors. In this section, we discuss the proper way of defining the corresponding realization of $\mbox{\normalfont SO}_0$ in $T(V)/I$, i.e., the analogue of the Spin group in this framework.

\subsection{Supervector reflections}
 The group of linear transformations generated by the supervector reflections was briefly introduced in \cite{sommen2001clifford} using the notion of the unit sphere in $\R^{m,2n}(\Lam_N)$ defined as $\mathbb S(m|2n)(\Lam_N)=\{{\bf w}\in \R^{m,2n}(\Lam_N): {\bf w}^2=-1\}.$
 The {\it reflection} associated to the supervector ${\bf w}\in \mathbb S(m|2n)(\Lam_N)$ is defined by the linear action on supervector variables  
\begin{equation}\label{Refl}
\psi({\bf w})[{\bf x}]={\bf w}{\bf x}{\bf w}, \hspace{2cm} {\bf x}\in{\bf S}.
\end{equation}
 It is known from the radial algebra setting that $\psi({\bf w})[{\bf x}]$ yields a new supervector variable. Indeed, for $x,y\in S$ one has ${ y}{ x}{ y}=\{{ x},{ y}\}{ y}-{ y}^2{ x}=\{{ x},{ y}\}{ y}+{ x}$.
 Every supervector reflection can be represented by a supermatrix in $\Mat(m|2n)(\Lam_N)$.
\begin{lem}\label{SupVecRefInsuf}
Let ${\bf w}=\underline{w}+\underline{w\p}=\sum_{j=1}^m w_je_j+ \sum_{j=1}^{2n} w\p_je\p_j\in\mathbb S(m|2n)(\Lam_N)$. Then, the linear transformation (\ref{Refl}) can be represented by a supermatrix 
\[\psi({\bf w})=\left(\begin{array}{cc} A({\bf w}) & B^{\p}({\bf w}) \\ C^{\p}({\bf w}) & D({\bf w})\end{array}\right)\in \Mat(m|2n)(\Lam_N)\]
with 
$A({\bf w}) = -2D_{\underline{w}}E_{m\times m}D_{\underline{w}}+I_m$, $B^{\p}({\bf w}) = D_{\underline{w}} E_{m\times 2n} D_{\underline{w\p}} \,\Om_{2n}$, $C^{\p}({\bf w}) = -2D_{\underline{w\p}}E_{2n\times m}D_{\underline{w}}$, and finally $D({\bf w}) = D_{\underline{w\p}}E_{2n\times 2n}D_{\underline{w\p}}\,\Om_{2n}+I_{2n}$, 
where 
{\footnotesize 
\begin{align*}
D_{\underline{w}}&=\left(\hspace{-.2cm}\begin{array}{ccc} w_1 & & \\  & \ddots & \\ & & w_m\end{array}\hspace{-.2cm}\right), & 
D_{\underline{w\p}}&=\left(\hspace{-.2cm} \begin{array}{ccc} w\p_1 & & \\  & \ddots & \\ & & w\p_{2n}\end{array}\hspace{-.2cm} \right), &
E_{p\times q}&=\left(\hspace{-.2cm} \begin{array}{cccc} 1 & 1& \ldots & 1 \\    \vdots & \vdots& \ddots & \vdots\\ 1 & 1& \ldots & 1\end{array}\hspace{-.2cm} \right)\in \R^{p\times q}.
\end{align*}
}
\end{lem}
\pf
Observe that $\psi({\bf w})[{\bf x}] = {\bf w}{\bf x}{\bf w}= \{{\bf x},{\bf w}\}{\bf w}+{\bf x}=\sum_{k=1}^m y_ke_k+ \sum_{k=1}^{2n} y\p_ke\p_k$, where
\begin{align*}
y_k &= -2\left(\sum_{j=1}^m w_jw_k x_j\right)+ x_k + \sum_{j=1}^n w\p_{2j-1}w_k x\p_{2j}-w\p_{2j}w_k x\p_{2j-1}, \\
y\p_k &= -2\left(\sum_{j=1}^m w_jw\p_k x_j\right)+ x\p_k + \sum_{j=1}^n -w\p_{2j-1}w\p_k x\p_{2j}+w\p_{2j}w\p_k x\p_{2j-1}.
\end{align*}
Then, $\psi({\bf w}) {\bf x}=\left(\begin{array}{cc} A({\bf w}) & B^{\p}({\bf w}) \\ C^{\p}({\bf w}) & D({\bf w})\end{array}\right) \left(\begin{array}{c} \underline{x} \\ \underline{x\p}\end{array}\right)$ 
where, 
{\footnotesize 
\begin{align*}
A(w) &= -2\left(\begin{array}{cccc} w_1^2 & w_2w_1& \ldots & w_mw_1 \\  w_1w_2 & w_2^2& \ldots & w_mw_2 \\  \vdots & \vdots& \ddots & \vdots\\ w_1w_m & w_2w_m& \ldots & w_m^2\end{array}\right)+I_m=-2D_{\underline{w}}E_{m\times m}D_{\underline{w}}+I_m,\\
B^{\p}(w) &= \left(\begin{array}{ccccc} -w\p_2w_1 & w\p_1w_1 & \ldots & -w\p_{2n}w_1 & w\p_{2n-1}w_1 \\  -w\p_2w_2 & w\p_1w_2 & \ldots & -w\p_{2n}w_2 & w\p_{2n-1}w_2  \\  \vdots & \vdots& \ddots & \vdots & \vdots \\ -w\p_2w_m & w\p_1w_m & \ldots & -w\p_{2n}w_m & w\p_{2n-1}w_m \end{array}\right)=D_{\underline{w}}E_{m\times 2n}D_{\underline{w\p}} \,J_{2n},\\
C^{\p}(w) &= -2\left(\begin{array}{cccc} w_1w\p_1 & w_2w\p_1& \ldots & w_mw\p_1 \\  w_1w\p_2 & w_2w\p_2& \ldots & w_mw\p_2 \\  \vdots & \vdots& \ddots & \vdots\\ w_1w\p_{2n} & w_2w\p_{2n}& \ldots & w_mw\p_{2n}\end{array}\right)=-2D_{\underline{w\p}}E_{2n\times m}D_{\underline{w}},\\
D(w) &= \left(\begin{array}{ccccc} w\p_2w\p_1 & -w\p_1w\p_1 & \ldots & w\p_{2n}w\p_1 & -w\p_{2n-1}w\p_1 \\  w\p_2w\p_2 & -w\p_1w\p_2 & \ldots & w\p_{2n}w\p_2 & -w\p_{2n-1}w\p_2  \\  \vdots & \vdots& \ddots & \vdots & \vdots \\ w\p_2w\p_{2n} & -w\p_1w\p_{2n} & \ldots & w\p_{2n}w\p_{2n} & -w\p_{2n-1}w\p_{2n} \end{array}\right)+I_{2n}=D_{\underline{w\p}}E_{2n\times 2n}D_{\underline{w\p}}\,J_{2n}+I_{2n}.
\end{align*}
}
$\hfill\square$

Algebraic operations with the matrices $A(w), B^{\p}(w), C^{\p}(w), D(w)$ are easy since 
\begin{equation}\label{Mat1}
E_{p\times m} D_{\underline{w}}^2 E_{m\times q} = \left(\sum_{j=1}^m w_j^2\right)E_{p\times q}, \hspace{.5cm} E_{p\times 2n} D_{\underline{w\p}}J_{2n}D_{\underline{w\p}} E_{2n\times q} = 2\left(\sum_{j=1}^n w\p_{2j-1}w\p_{2j}\right)E_{p\times q}.
\end{equation}

\begin{pro}
Let ${\bf w}\in \mathbb S(m|2n)(\Lam_N)$. Then $\psi({\bf w})\in \mbox{\normalfont O}_0$ and $\sd\left(\psi({\bf w})\right)=-1$.
\end{pro}
\pf
In order to prove that $\psi({\bf w})\in \mbox{\normalfont O}_0$ it suffices to prove that $A({\bf w}), B^{\p}({\bf w}), C^{\p}({\bf w}), D({\bf w})$ satisfy (\ref{ConCmp}). This can be easily done using (\ref{Mat1}) and the identity $-1={\bf w}^2=-\sum_{j=1}^m w_j^2+\sum_{j=1}^n w\p_{2j-1}w\p_{2j}$. In fact, we have 
\begin{align*}
A({\bf w})^TA({\bf w}) &= 4D_{\underline{w}}E_{m\times m}D_{\underline{w}}^2E_{m\times m}D_{\underline{w}}-4D_{\underline{w}}E_{m\times m}D_{\underline{w}}+I_m\\
&=4\left(\sum_{j=1}^m w_j^2\right)D_{\underline{w}}E_{m\times m}D_{\underline{w}}-4D_{\underline{w}}E_{m\times m}D_{\underline{w}}+I_m,
\end{align*} 
and $C^{\p}({\bf w})^T\Om_{2n}C^{\p}({\bf w}) = 4D_{\underline{w}}E_{m\times 2n}D_{\underline{w\p}}\Om_{2n}D_{\underline{w\p}}E_{2n\times m}D_{\underline{w}}=8\left(\sum_{j=1}^n w\p_{2j-1}w\p_{2j}\right)D_{\underline{w}}E_{m\times m}D_{\underline{w}}.$
Then, $A({\bf w})^TA({\bf w})-\frac{1}{2}C^{\p}({\bf w})^T\Om_{2n}C^{\p}({\bf w})=4\left[\sum_{j=1}^m w_j^2-1-\sum_{j=1}^n w\p_{2j-1}w\p_{2j}\right]D_{\underline{w}}E_{m\times m}D_{\underline{w}}+I_m=I_m$.
Also,
\begin{align*}
A({\bf w})^TB^{\p}({\bf w}) &= -2D_{\underline{w}}E_{m\times m}D_{\underline{w}}^2E_{m\times 2n}D_{\underline{w\p}}\Om_{2n}+D_{\underline{w}}E_{m\times 2n}D_{\underline{w\p}} \,\Om_{2n}\\
&=-2\left(\sum_{j=1}^m w_j^2\right)D_{\underline{w}}E_{m\times 2n}D_{\underline{w\p}}\Om_{2n}+D_{\underline{w}}E_{m\times 2n}D_{\underline{w\p}} \,\Om_{2n},
\end{align*} 
and
\begin{align*}
C^{\p}({\bf w})^T\Om_{2n}D({\bf w}) &= -2D_{\underline{w}}E_{m\times 2n}D_{\underline{w\p}}\Om_{2n}D_{\underline{w\p}}E_{2n\times 2n}D_{\underline{w\p}}\Om_{2n}-2D_{\underline{w}}E_{m\times 2n}D_{\underline{w\p}}\Om_{2n}\\
&=-4\left(\sum_{j=1}^n w\p_{2j-1}w\p_{2j}\right)D_{\underline{w}}E_{m\times 2n}D_{\underline{w\p}}\Om_{2n}-2D_{\underline{w}}E_{m\times 2n}D_{\underline{w\p}}\Om_{2n}.
\end{align*} 
Hence, $A({\bf w})^TB^{\p}({\bf w})-\frac{1}{2}C^{\p}({\bf w})^T\Om_{2n}D({\bf w}) = 2\left[-\sum_{j=1}^m w_j^2+1+\sum_{j=1}^n w\p_{2j-1}w\p_{2j}\right]D_{\underline{w}}E_{m\times 2n}D_{\underline{w\p}}\Om_{2n}=0$.
 In the same way we have
 \begin{align*}
B^{\p}({\bf w})^TB^{\p}({\bf w}) &= -\Om_{2n}D_{\underline{w\p}}E_{2n\times m}D_{\underline{w}}^2E_{m\times 2n}D_{\underline{w\p}}\Om_{2n}=-\left(\sum_{j=1}^m w_j^2\right)\Om_{2n}D_{\underline{w\p}}E_{2n\times 2n}D_{\underline{w\p}}\Om_{2n},
\end{align*} 
and
\begin{align*}
D({\bf w})^T\Om_{2n}D({\bf w}) &= \Om_{2n}D_{\underline{w\p}}E_{2n\times 2n}D_{\underline{w\p}}\Om_{2n}D_{\underline{w\p}}E_{2n\times 2n}D_{\underline{w\p}}\Om_{2n}+2\Om_{2n}D_{\underline{w\p}}E_{2n\times 2n}D_{\underline{w\p}}\Om_{2n}+\Om_{2n}\\
&=2\left(\sum_{j=1}^n w\p_{2j-1}w\p_{2j}+1\right)\Om_{2n}D_{\underline{w\p}}E_{2n\times 2n}D_{\underline{w\p}}\Om_{2n}+\Om_{2n},
\end{align*} 
 whence
\begin{align*}
B^{\p}({\bf w})^TB^{\p}({\bf w})\hspace{-.1cm} + \hspace{-.1cm} \frac{1}{2}D^{T}({\bf w})\Om_{2n}D({\bf w}) \hspace{-.1cm}&= \hspace{-.1cm} \left[ \hspace{-.1cm}-\sum_{j=1}^m w_j^2+\hspace{-.1cm}1\hspace{-.1cm}+\sum_{j=1}^n w\p_{2j-1}w\p_{2j}\hspace{-.1cm} \right]\hspace{-.1cm} \Om_{2n}D_{\underline{w\p}}E_{2n\times 2n}D_{\underline{w\p}}\Om_{2n}+\frac{1}{2}\Om_{2n} = \frac{1}{2}\Om_{2n}.
\end{align*}
 Then, $A({\bf w}), B^{\p}({\bf w}), C^{\p}({\bf w}), D({\bf w})$ satisfy (\ref{ConCmp}) and in consequence, $\psi({\bf w})\in \mbox{\normalfont O}_0$.  To prove that $\sd(\psi({\bf w}))=-1$, first observe that $\psi({\bf w})=\psi({\bf w})^{-1}$ since $\psi({\bf w})\circ \psi({\bf w})[{\bf x}]={\bf w}{\bf w}{\bf x}{\bf w}{\bf w}={\bf x}$.
 Hence, due to (\ref{Inv}) we obtain
$A({\bf w})=\left(A({\bf w})-B^{\p}({\bf w})D({\bf w})^{-1}C^{\p}({\bf w})\right)^{-1}, $
 yielding
 \[\sd(\psi({\bf w}))=\frac{\det\left[A({\bf w})-B^{\p}({\bf w})D({\bf w})^{-1}C^{\p}({\bf w})\right]}{\det[D({\bf w})]}=\frac{1}{\det[A({\bf w})]\det[D({\bf w})]}.\]
 We will compute $\det[D({\bf w})]$ using the formula $\det[D({\bf w})]=\exp(\tr \ln D({\bf w}))$ and the fact that $D({\bf w})-I_{2n}$ is a nilpotent matrix. Observe that 
\[\ln D({\bf w})=\sum_{j=1}^\infty(-1)^{j+1}\frac{\left(D({\bf w})-I_{2n}\right)^j}{j}=\sum_{j=1}^\infty(-1)^{j+1}\frac{\left(D_{\underline{w\p}}E_{2n\times 2n}D_{\underline{w\p}}\,\Om_{2n}\right)^j}{j}.\]
It follows from (\ref{Mat1}) that  $\left(D_{\underline{w\p}}E_{2n\times 2n}D_{\underline{w\p}}\,\Om_{2n}\right)^j=2^{j-1}\left(\sum_{j=1}^n w\p_{2j-1}w\p_{2j}\right)^{j-1}D_{\underline{w\p}}E_{2n\times 2n}D_{\underline{w\p}}\,\Om_{2n}$.
Then,
\[\ln D({\bf w})= \left[\sum_{j=1}^\infty(-1)^{j+1}\frac{2^{j-1}\left(\sum_{j=1}^n w\p_{2j-1}w\p_{2j}\right)^{j-1}}{j}\right]D_{\underline{w\p}}E_{2n\times 2n}D_{\underline{w\p}}\,\Om_{2n}\]
and in consequence, 
\[\displaystyle \tr \ln D({\bf w})=-\sum_{j=1}^\infty(-1)^{j+1}\frac{2^{j}\left(\sum_{j=1}^n w\p_{2j-1}w\p_{2j}\right)^{j}}{j}=-\ln\left(1+2\sum_{j=1}^n w\p_{2j-1}w\p_{2j}\right).\]
Hence $\det(D({\bf w}))=\left(1+2\sum_{j=1}^n w\p_{2j-1}w\p_{2j}\right)^{-1}$.
Similar computations yield $\det(A({\bf w}))=1-2\sum\limits_{j=1}^m w_j^2$, 
which finally shows that $\sd(\psi({\bf w}))=-1$.
$\hfill\square$
 
We can now define the {\it bosonic Pin group} in this setting as
\[\Pin_b(m|2n)(\Lam_N) = \{{\bf w}_1\cdots {\bf w}_k \,:\, {\bf w}_j\in \mathbb S(m|2n)(\Lam_N), k\in \N\},\]
and extend the map $\psi$ to a Lie group homomorphism $\psi: \Pin_b(m|2n)(\Lam_N)\fd \mbox{\normalfont O}_0$ by
\[\psi({\bf w}_1\cdots {\bf w}_k)[{\bf x}]={\bf w}_1\cdots {\bf w}_k \,{\bf x}\, {\bf w}_k\cdots {\bf w}_1=\psi({\bf w}_1)\circ\cdots\circ\psi({\bf w}_k)[{\bf x}].\]
It is clearly seen that the restriction of $\psi$ to the {\it bosonic spin group}, defined as
\[\Spin_b(m|2n)(\Lam_N)=\{{\bf w}_1\cdots {\bf w}_{2k}:{\bf w}_j\in \mathbb S(m|2n)(\Lam_N), k\in \N\},\]
takes values in the subgroup $\mbox{\normalfont SO}_0\inc \mbox{\normalfont O}_0$.
 
In the classical case, the Pin group and the Spin group are double coverings of the groups $\mbox{\normalfont O}(m)$ and $\mbox{\normalfont SO}(m)$ respectively. A natural question in this setting is whether $\Pin_b(m|2n)(\Lam_N)$ and $\Spin_b(m|2n)(\Lam_N)$ cover the groups $\mbox{\normalfont O}_0$ and $\mbox{\normalfont SO}_0$.  The answer to this question is negative and the main reason for this is that the real projection of every vector ${\bf w}\in \Sa(m|2n)(\Lam_N)$ is in the unitary sphere $\Sa^{m-1}$ of $\R^m$, i.e.,
\[[{\bf w}]_0=\sum_{j=1}^m[w_j]_0 e_j \hspace{.5cm}\mbox{ and }\hspace{.5cm} [{\bf w}]_0^2=-1.\] 
Then, the real projection of $\psi\left(\Pin_b(m|2n)(\Lam_N)\right)$ is just $\mbox{\normalfont O}(m)$, while $[\mbox{\normalfont O}_0]_0=\mbox{\normalfont O}(m)\times \mbox{\normalfont Sp}_\Om(2n)$. This means that these bosonic versions of Pin and Spin do not describe the symplectic parts of $\mbox{\normalfont O}_0$ and $\mbox{\normalfont SO}_0$. This phenomenon is due to the natural structure of supervectors: their real projections belong to a space with an orthogonal structure while the symplectic structure plays no r\^ole. Up to a nilpotent vector, they are classical Clifford vectors, whence it is impossible to generate by this approach the real symplectic geometry that is also present in the structure of $\mbox{\normalfont O}_0$ and $\mbox{\normalfont SO}_0$. That is why we have chosen the name of "bosonic" Pin and "bosonic" Spin groups. This also explains why we had to extend the space of superbivectors in Section \ref{SupBi-VecExtSect}. The ordinary superbivectors in $\Lam_N\otimes \mathcal{C}_{m,2n}$ are generated over $\Lam_N^{(ev)}$ by the wedge product of supervectors. Then, they can only describe $\mathfrak{so}(m)$ and not $\mathfrak{sp}_\Om(2n)$ and in consequence, they do not cover $\mathfrak{so}_0$.
 
 As in the classical setting (see \cite{MR1777332}), it is possible to obtain the following result that shows, from another point of view, that $\Pin_b(m|2n)(\Lam_N)$ cannot completely describe $\mbox{\normalfont O}_0$.
 \begin{pro}
The Lie algebra of $\Pin_b(m|2n)(\Lam_N)$ is included in $\R_{m|2n}^{(2)}(\Lam_N)$.
\end{pro}
\pf
Let $\gam(t)={\bf w}_1(t)\cdots {\bf w}_k(t)$ be a path in $\Pin_b(m|2n)(\Lam_N)$ with ${\bf w}_j(t)\in  \mathbb S(m|2n)(\Lam_N)$ for every $t\in \R$ and $\gam(0)=1$.  The tangent to $\gam$ at $t=0$ is $\frac{d\gam}{dt}\big|_{t=0}= \sum_{j=1}^k {\bf w}_1(0)\cdots {\bf w}'_j(0)\cdots {\bf w}_k(0)$. We will show that each summand of $\frac{d\gam}{dt}\big|_{t=0}$ belongs to $\R_{m|2n}^{(2)}(\Lam_N)$.

\noindent For $j=1$ we have ${\bf w}'_1(0) {\bf w}_2(0)\cdots {\bf w}_k(0)=-{\bf w}'_1(0){\bf w}_1(0)$. But ${\bf w}_1(t){\bf w}_1(t)\equiv -1$ implies
\[\{{\bf w}'_1(0),{\bf w}_1(0)\}={\bf w}'_1(0){\bf w}_1(0)+{\bf w}_1(0){\bf w}'_1(0)=0. 
\]
Then ${\bf w}'_1(0){\bf w}_1(0)=\frac{1}{2}\{{\bf w}'_1(0),{\bf w}_1(0)\}+{\bf w}'_1(0)\wedge {\bf w}_1(0)={\bf w}'_1(0)\wedge {\bf w}_1(0)\in \R_{m|2n}^{(2)}(\Lam_N)$. For $j=2$, 
\begin{align*}
{\bf w}_1(0) {\bf w}'_2(0)\cdots {\bf w}_k(0) &={\bf w}_1(0) {\bf w}'_2(0){\bf w}_2(0) {\bf w}_1(0) = \ -\left[{\bf w}_1(0) {\bf w}'_2(0){\bf w}_1(0)\right] \, \left[{\bf w}_1(0) {\bf w}_2(0){\bf w}_1(0)\right]\\
&=-\psi({\bf w}_1(0))[{\bf w}'_2(0)]\; \psi({\bf w}_1(0))[{\bf w}_2(0)].
\end{align*}
But $\psi({\bf w}_1(0))\in \mbox{\normalfont O}_0$ preserves the inner product (see remark \ref{InvGrasInProdRem}), so 
\[{\bf w}_1(0) {\bf w}'_2(0)\cdots {\bf w}_k(0) = \psi({\bf w}_1(0))[{\bf w}'_2(0)]\wedge \psi({\bf w}_1(0))[{\bf w}_2(0)]\in \R_{m|2n}^{(2)}(\Lam_N).\]
We can proceed similarly for every $j=3,\ldots, k$.
$\hfill\square$
 
\subsection{A proper definition for the group $\Spin(m|2n)(\Lam_N)$ }
The above approach shows that the radial algebra setting does not contain a suitable realization of $\mbox{\normalfont SO}_0$ in the Clifford superspace framework.  Observe that the Clifford representation of $\mathfrak{so}_0$ given by $\R_{m|2n}^{(2)E}(\Lam_N)$ lies outside of the radial algebra $\R_{m|2n}(\Lam_N)$, which suggests that something similar should happen with 
the corresponding Lie group $\mbox{\normalfont SO}_0$. In this case, a proper definition for the Spin group would be generated by the exponentials (in general contained in $T(V)/I$) of all elements in $\R_{m|2n}^{(2)E}(\Lam_N)$, i.e.
\[\Spin(m|2n)(\Lam_N):=\left\{e^{B_1}\cdots e^{B_k}: B_1, \ldots, B_k\in \R_{m|2n}^{(2)E}(\Lam_N), k\in\N\right\},\]
and the action of this group on supervector variables ${\bf x}\in{\bf S}$ 
 is given by the group homomorphism $h: \Spin(m|2n)(\Lam_N)\fd \mbox{\normalfont SO}_0$ defined by  
 \begin{equation}\label{Rep_h_SpinGSS}
 h(e^B)[{\bf x}]=e^B {\bf x} e^{-B}, \hspace{1cm} B\in \R_{m|2n}^{(2)E}(\Lam_N), \; {\bf x}\in {\bf S}.
 \end{equation}
 In fact, for every extended superbivector $B$, $h(e^B)$ maps  supervector variables into new supervector variables and admits a supermatrix representation in $\Mat(m|2n)(\Lam_N)$ belonging to $\mbox{\normalfont SO}_0$. This is summarized below.
  \begin{pro}\label{InfRep}
Let $B\in \R_{m|2n}^{(2)E}(\Lam_N)$. Then, 
$h(e^B)[{\bf x}]=e^{\phi(B)}{\bf x}$ for every $ {\bf x}\in {\bf S}$.
\end{pro}
\pf
In the associative algebra $\mathcal{A}_{m,2n}\otimes \Lam_N$, the identity 
$\displaystyle \underbrace{[B,[B\ldots[B,{\bf x}]\ldots]]}_k=\sum_{j=0}^k \binom{k}{j}B^j{\bf x}(-B)^{k-j}$
holds. Then,
\begin{align*}
h(e^B)[{\bf x}] &= e^B {\bf x} e^{-B} \ = \ \left(\sum_{k=0}^\infty \frac{B^k}{k!}\right){\bf x} \left(\sum_{k=0}^\infty \frac{(-B)^k}{k!}\right)  = \sum_{k=0}^\infty \left(\sum_{j=0}^k \frac{B^j}{j!}\,{\bf x}\,\frac{(-B)^{k-j}}{(k-j)!}\right)\\
&= \sum_{k=0}^\infty \frac{1}{k!}\left(\sum_{j=0}^k\binom{k}{j}B^j{\bf x}(-B)^{k-j}\right)  =\sum_{k=0}^\infty \frac{1}{k!}\underbrace{[B,[B\ldots[B,{\bf x}]\ldots]]}_k = \sum_{k=0}^\infty \frac{\phi(B)^k {\bf x}}{k!} = e^{\phi(B)}{\bf x}.
\end{align*}
$\hfill\square$
\begin{remark}
Proposition \ref{InfRep} means that the Lie algebra isomorphism $\phi: \R_{m|2n}^{(2)E}(\Lam_N)\fd \mathfrak{so}_0$ is the derivative at the origin (or infinitesimal representation) of the Lie group homomorphism $h: \Spin(m|2n)(\Lam_N)\fd \mbox{\normalfont SO}_0$, i.e.,
\begin{equation}\label{psiAsInfRepofH}
e^{t\phi(B)}=h(e^{tB}) \hspace{1.cm} \pt t\in\R, \;  B\in \R_{m|2n}^{(2)E}(\Lam_N).
\end{equation}
\end{remark}
On account of the connectedness of $\mbox{\normalfont SO}_0$ it can be shown  that the group $\Spin(m|2n)(\Lam_N)$ is a realization of $\mbox{\normalfont SO}_0$ in $T(V)/I$ through the representation $h$.
\begin{teo}
For every $M\in \mbox{\normalfont SO}_0$ there exists an element $s\in \Spin(m|2n)(\Lam_N)$ such that $h(s)=M$. 
\end{teo}
\pf Since $\mbox{\normalfont SO}_0$ is a connected Lie group (Proposition \ref{Con}), for every supermatrix $M\in \mbox{\normalfont SO}_0$ there exist $X_1,\ldots, X_k\in \mathfrak{so}_0$ such that $e^{X_1}\cdots e^{X_k}=M$, see Corollary 3.47 in \cite{MR3331229}. Taking $B_1, \ldots, B_k \in \R_{m|2n}^{(2)E}(\Lam_N)$ such that $\phi(B_j)=X_j$, $j=1,\ldots, k$, we obtain
\[M{\bf x} = e^{X_1}\cdots e^{X_k} \, {\bf x}= e^{\phi(B_1)}\cdots e^{\phi(B_k)}\, {\bf x} =  h(e^{B_1})\circ \cdots \circ h(e^{B_k})[{\bf x}] = h( e^{B_1}\cdots e^{B_k})[{\bf x}].\]
Then, 
 $s= e^{B_1}\cdots e^{B_k}\in \Spin(m|2n)(\Lam_N)$ satisfies $h(s)=M$.
$\hfill\square$

The decomposition of $\mbox{\normalfont SO}_0$ given in Theorem \ref{SO_0 Decomp} provides the exact number of exponentials of extended superbivectors to be considered in $\Spin(m|2n)(\Lam_N)$ in order to cover the whole group $\mbox{\normalfont SO}_0$. 
If we consider the subspaces $\Xi_1, \Xi_2, \Xi_3$ of $\R_{m|2n}^{(2)E}(\Lam_N)$ given by
\begin{align}
\Xi_1 &= \phi^{-1}\big(\mathfrak{so}(m)\times [ \mathfrak{sp}_\Om(2n)\cap \mathfrak{so}(2n)]\big), & \dim \Xi_1&=\frac{m(m-1)}{2}+n^2,\nonumber\\
\Xi_2 &= \phi^{-1}\big(\{0_m\}\times[\mathfrak{sp}_\Om(2n)\cap \Sym(2n)]\big),                                & \dim \Xi_2&=n^2+n,\label{SpinDec}\\
\Xi_3 &= \phi^{-1}\big(\mathfrak{so}_0(m|2n)(\Lam_N^+)\big)                 & \dim \Xi_3&=\dim \mathfrak{so}_0 - \frac{m(m-1)}{2} - n(2n+1), \nonumber
\end{align}
we get the decomposition $\R_{m|2n}^{(2)E}(\Lam_N)=\Xi_1\oplus \Xi_2\oplus \Xi_3$, leading to the subset  \[\Xi=\exp(\Xi_1) \exp(\Xi_2) \exp(\Xi_3)\inc \Spin(m|2n)(\Lam_N),\] which suffices for describing $\mbox{\normalfont SO}_0$.  Indeed, from Theorem \ref{SO_0 Decomp} it follows that the restriction $h: \Xi\fd \mbox{\normalfont SO}_0$ is surjective. 
We now investigate the explicit form of the superbivectors in each of the subspaces $\Xi_1$, $\Xi_2$ and $\Xi_3$. 
\begin{pro}\label{S_jBasis} 
The following statements hold.
\begin{itemize}
\item[(i)] A basis for $\Xi_1$ is $\begin{cases} e_je_k,  &1\leq j< k\leq m,\\ e\p_{2j-1}\odot e\p_{2k-1}+e\p_{2j}\odot e\p_{2k}, & 1\leq j\leq k\leq n,\\ e\p_{2j-1}\odot e\p_{2k}-e\p_{2j}\odot e\p_{2k-1},  &1\leq j< k\leq n. \end{cases}$
\item[(ii)] A basis for $\Xi_2$ is: $\begin{cases} e\p_{2j-1}\odot e\p_{2j},  &1\leq j\leq n,\\ e\p_{2j-1}\odot e\p_{2k-1}-e\p_{2j}\odot e\p_{2k}, & 1\leq j\leq k\leq n,\\ e\p_{2j-1}\odot e\p_{2k}+e\p_{2j}\odot e\p_{2k-1},  &1\leq j< k\leq n. \end{cases}$
\item[(iii)] $\Xi_3$ consists of all elements of the form (\ref{SBiv}) with $b_{j,k}, B_{j,k}\in \Lam_N^{(ev)}\cap \Lam_N^{+}$ and  $b\p_{j,k}\in \Lam_N^{(odd)}$.
\end{itemize}
\end{pro}
\pf We first recall that a basis for the Lie algebra $\mathfrak{sp}_\Om(2n)$ is given by the elements
\begin{align*}
A_{j,k} &:=E_{2j,2k-1}+E_{2k,2j-1}, \hspace{.2cm}1\leq j\leq k\leq n & B_{j,k}&:=E_{2j-1,2k}+E_{2k-1,2j}, \hspace{.2cm}1\leq j\leq k\leq n,\\
C_{j,k} &:=E_{2k,2j}-E_{2j-1,2k-1}, \hspace{.2cm}1\leq j\leq k\leq n, & D_{j,k} &:=E_{2j,2k}-E_{2k-1,2j-1}, \hspace{.2cm} 1\leq j<k\leq n,
\end{align*} 
where the matrices $E_{j,k}\in \R^{n\times n}$ are defined as in Lemma \ref{LemIso}. It holds that $A_{j,k}^T=B_{j,k}$ for $1\leq j\leq k\leq n$, $C_{j,k}^T=D_{j,k}$ for $1\leq j< k\leq n$ and $C_{j,j}^T=C_{j,j}$ for $1\leq j\leq n$. Hence, for every $D_0\in \mathfrak{sp}_\Om(2n)$ we have
\begin{align*}
D_0 &= \sum_{1\leq j\leq k\leq n} \left(a_{j,k}A_{j,k}+b_{j,k}B_{j,k}+ c_{j,k}C_{j,k}\right) + \sum_{1\leq j< k\leq n} d_{j,k}D_{j,k},\\
D_0^T &= \sum_{1\leq j\leq k\leq n} \left(a_{j,k}B_{j,k}+b_{j,k}A_{j,k}\right) + \sum_{1\leq j< k\leq n} \left(c_{j,k}D_{j,k} + d_{j,k}C_{j,k}\right) + \sum_{j=1}^n c_{j,j}C_{j,j},
\end{align*}
where $a_{j,k},b_{j,k}, c_{j,k}, d_{j,k}\in \R$.
\begin{itemize}
\item[$(i)$] From the previous equalities we get that $D_0^T=-D_0$ if and only if  
 \[\displaystyle D_0=\sum_{1\leq j\leq k\leq n} a_{j,k}\left(A_{j,k}-B_{j,k}\right) + \sum_{1\leq j< k\leq n} c_{j,k}\left(C_{j,k} -D_{j,k}\right).\]
  Then, a basis for $ \mathfrak{sp}_\Om(2n)\cap \mathfrak{so}(2n)$ is $\{A_{j,k}-B_{j,k}: 1\leq j\leq k\leq n\}\cup \{C_{j,k} -D_{j,k}:1\leq j< k\leq n\}$.
The remainder of the proof directly follows from Lemma \ref{LemIso}.
\item[$(ii)$] In this case we have that $D_0^T=D_0$ if and only if 
\[\displaystyle D_0=\sum_{1\leq j\leq k\leq n} a_{j,k}\left(A_{j,k}+B_{j,k}\right) + \sum_{1\leq j< k\leq n} c_{j,k}\left(C_{j,k} +D_{j,k}\right)+ \sum_{j=1}^n c_{j,j} C_{j,j},\]
whence a basis for $\mathfrak{sp}_\Om(2n)\cap \Sym(2n)$ is
\[
\{A_{j,k}+B_{j,k}: 1\leq j\leq k\leq n\}\cup \{C_{j,j}:1\leq j\leq n\} \cup \{C_{j,k} +D_{j,k}:1\leq j< k\leq n\}.\]
The remainder of the proof directly follows from Lemma \ref{LemIso}.
\item[iii)] This trivially follows from Lemma \ref{LemIso}. $\hfill\square$
\end{itemize}

\subsection{Spin covering of the group $\mbox{\normalfont SO}_0$}
It is a natural question in this setting whether the spin group still is a double covering of the group of rotations, as it is in classical Clifford analysis. In other words, we will investigate how many times $\Xi\inc \Spin(m|2n)(\Lam_N)$ covers $\mbox{\normalfont SO}_0$, or more precisely, we will determine the cardinality of the set $\{s\in \Xi: h(s)=M\}$ given a certain fixed element $M\in \mbox{\normalfont SO}_0$.

From Proposition \ref{InfRep} we have that the representation $h$ of an element $s=e^{B_1} e^{B_2} e^{B_3}\in \Xi$, $B_j\in \Xi_j$, 
has the form $h(s)=e^{\phi(B_1)}e^{\phi(B_2)}e^{\phi(B_3)}$. 
Following the decomposition $M=e^Xe^Ye^{\bf Z}$ 
given in Theorem \ref{SO_0 Decomp} for $M\in \mbox{\normalfont SO}_0$, we get that $h(s)=M$ if and only if $e^{\phi(B_1)}=e^X$, $B_2=\phi^{-1}(Y)$ and $B_3=\phi^{-1}({\bf Z})$.
Then, the cardinality of $\{s\in \Xi: h(s)=M\}$ only depends on the number of extended superbivectors $B_1\in \Xi_1$ that satisfy $e^{\phi(B_1)}=e^X$.  It reduces our analysis to finding the kernel of the restriction  $h|_{\exp(\Xi_1)}:\exp(\Xi_1)\fd \mbox{\normalfont SO}(m)\times [\mbox{\normalfont Sp}_\Om(2n)\cap \mbox{\normalfont SO}(2n)]$ 
 of the Lie group homomorphism $h$ to $\exp(\Xi_1)$.  This kernel is given by $\ker h|_{\exp(\Xi_1)}=\{e^{B}\;\;:\;\; e^{\phi(B)}=I_{m+2n}, \,B\in \Xi_1\}$.
 
  We recall, from Proposition \ref{S_jBasis}, that $B\in \Xi_1$ may be written as $B=B_o+B_s$ where $B_o\in \R_{0,m}^{(2)}$ is a classical real bivector and $B_s\in \phi^{-1}\left(\{0_m\}\times [ \mathfrak{sp}_\Om(2n)\cap \mathfrak{so}(2n)]\right)$. 
 The components $B_o, B_s$ commute and in consequence, $e^B=e^{B_o} e^{B_s}$. Consider the projections $\phi_o$ and $\phi_s$ of $\phi$ over the algebra of classical bivectors $\R_{0,m}^{(2)}$ and over the algebra $\phi^{-1}\left(\{0_m\}\times [ \mathfrak{sp}_\Om(2n)\cap \mathfrak{so}(2n)]\right)$ respectively, i.e.
 \[\phi_o: \R_{0,m}^{(2)}\fd \mathfrak{so}(m), \hspace{.5cm} \phi_s: \phi^{-1}\left(\{0_m\}\times [ \mathfrak{sp}_\Om(2n)\cap \mathfrak{so}(2n)]\right) \fd \mathfrak{sp}_\Om(2n)\cap \mathfrak{so}(2n),\]
where $\phi(B)=\left(\begin{array}{cc} \phi_o(B_o) & 0 \\ 0 & \phi_s(B_s)\end{array}\right)$ for $B\in \Xi_1$.
Or equivalently,
\begin{align*}
\begin{cases}\phi_o(B_o)[\underline{x}]=[B_o,\underline{x}], \\ \phi_s(B_s)[\underline{x}\p]=\left[B_s,\underline{x}\p\right], \end{cases} && {\bf x}&=\underline{x}+\underline{x}\p\in {\bf S}.
\end{align*}
 Hence $e^{\phi(B)}=I_{m+2n}$ if and only if $e^{\phi_o(B_o)}=I_m$ and $e^{\phi_s(B_s)}=I_{2n}$. For the first condition, we know from classical Clifford analysis that $\Spin(m)= \{e^{B}:B\in \R_{0,m}^{(2)}\}$ is a double covering of $\mbox{\normalfont SO}(m)$ and in consequence $e^{\phi_o(B_o)}=I_m$ 
implies $e^{B_0}=\pm 1$.
 Let us now compute all possible values for $e^{B_s}$ for which $e^{\phi_s(B_s)}=I_{2n}$. To that end, we need the following linear algebra result.
\begin{pro}
Every matrix $D_0\in \mathfrak{so}(2n)\cap \mathfrak{sp}_\Om(2n)$ can be written in the form $D_0=R\Sig R^T$ where $R\in \mbox{\normalfont SO}(2n)\cap \mbox{\normalfont Sp}_\Om(2n)$ and 
{
\begin{equation}\label{BlockDiag_SS}
\Sig = \left(\begin{array}{ccc} \begin{array}{cc} 0 & \cit_1 \\ -\cit_1 & 0\end{array} &  & \\  & \ddots & \\ & &  \begin{array}{cc} 0 & \cit_n \\ -\cit_n & 0\end{array} \end{array}\right), \hspace{1cm} \cit_j\in \R, \hspace{.3cm} j=1,\ldots, n. 
\end{equation}
}
\end{pro}
\pf
The map $\Psi(D_0)=\frac{1}{2} { Q}D_0 \left(Q^T\right)^c$, where
\[Q=\left(\begin{array}{ccccccc} 1 & i & 0 & 0 & \ldots & 0 & 0 \\ 0 & 0 & 1 & i & \ldots & 0 & 0 \\ \vdots & \vdots & \vdots & \vdots & \ddots & \vdots & \vdots  \\ 0 & 0 & 0 & 0 & \ldots & 1 & i \end{array}\right)\in \C^{n\times 2n},\]
is a Lie group isomorphism between $\mbox{\normalfont SO}(2n)\cap \mbox{\normalfont Sp}_\Om(2n)$ and $\mbox{U}(n)$. It is easily proven that $\Psi$ is its own infinitesimal representation on the Lie algebra level, and in consequence, a Lie algebra isomorphism between $ \mathfrak{so}(2n)\cap \mathfrak{sp}_\Om(2n)$ and $ \mathfrak{u}(n)$. The inverse of $\Psi$ is given by $\Psi^{-1}(L)=\frac{1}{2}\left[\left(Q^T\right)^c \,L \,Q + Q^T \,L^c \,Q^c\right]$.
 For every $D_0\in \mathfrak{so}(2n)\cap \mathfrak{sp}_\Om(2n)$, let us consider the skew-Hermitian matrix $L=\Psi(D_0)\in \mathfrak{u}(n)$.  It is  known that every skew-Hermitian matrix is unitarily diagonalizable and all its eigenvalues are purely imaginary, see \cite{MR2978290}.   Hence, $L=\Psi(D_0)$ can be written as $L=U\mathfrak{D} \left(U^T\right)^c$ 
 where $U\in \mbox{U}(n)$ and $\mathfrak{D}=\mbox{diag}(-i\cit_1,\ldots, -i\cit_n)$, $\cit_j\in \R$. 
 Then, $D_0=\Psi^{-1}(L)=R \Sig R^T$ 
 where $R=\Psi^{-1}(U)\in \mbox{\normalfont SO}(2n)\cap \mbox{\normalfont Sp}_\Om(2n)$ and $\Sig=\Psi^{-1}(\mathfrak{D})$ has the form  (\ref{BlockDiag_SS}).
$\hfill\square$
 
 Since $\phi_s(B_s)\in \mathfrak{so}(2n)\cap \mathfrak{sp}_\Om(2n)$, we have $\phi_s(B_s)= R \Sig R^T$ as in the previous proposition. Hence, $e^{\phi_s(B_s)}= R e^\Sig R^T$ where $e^\Sig$ is the block-diagonal matrix 
\[e^\Sig=\mbox{diag}(e^{\cit_1 \Om_2}, \ldots, e^{\cit_n \Om_2}) \quad \mbox{with\ } \quad e^{\cit_j \Om_2}=\cos \cit_j I_2+ \sin \cit_j \Om_2.\]
Hence $e^{\phi_s(B_s)}= I_{2n}$ if and only if $e^\Sig=I_{2n}$, which is seen to be equivalent to $\cit_j=2k_j\pi$,  $k_j\in\Z$ ($j=1,\ldots, n$), or to 
\[\Sig=\sum_{j=1}^n 2k_j\pi \left(E_{2j-1,2j}-E_{2j,2j-1}\right), \quad k_j\in\Z \; (j=1,\ldots, n).\]
Now, $\mbox{\normalfont SO}(2n)\cap \mbox{\normalfont Sp}_\Om(2n)$ being connected and compact, there exists $B_R\in \phi^{-1}(\mathfrak{so}(2n)\cap \mathfrak{sp}_\Om(2n))$ such that $R=e^{\phi(B_R)}$.  We recall that the $h$-action leaves any multivector structure invariant, in particular, $h[e^B]\left(\R_{m|2n}^{(2)E}(\Lam_N)\right)\inc \R_{m|2n}^{(2)E}(\Lam_N)$ 
for every $B\in\R_{m|2n}^{(2)E}(\Lam_N)$.  Then, using the fact that $\phi$ is the derivative at the origin of $h$, we get that the extended superbivector $h(e^{B_R})[\phi^{-1}(\Sig)]= e^{B_R}\phi^{-1}(\Sig)e^{-B_R}$ 
is such that
\[
\phi\left(e^{B_R}\phi^{-1}(\Sig)e^{-B_R}\right)=e^{\phi(B_R)}\Sig e^{-\phi(B_R)}=R\Sig R^T=\phi(B_s),
\]
implying that $B_s=e^{B_R}\phi^{-1}(\Sig)e^{-B_R}$.  Then, in order to compute $e^{B_s}=e^{B_R}\,e^{\phi^{-1}(\Sig)}\,e^{-B_R}$, we first have to compute $e^{\phi^{-1}(\Sig)}$. Following the correspondences given in Lemma \ref{LemIso} we get
\[\phi^{-1}(\Sig) = \sum_{j=1}^n 2k_j\pi \,\phi^{-1}\left(E_{2j-1,2j}-E_{2j,2j-1}\right)=  \sum_{j=1}^n k_j\pi \left(e\p_{2j-1}^{\,2}+e\p_{2j}^{\,2}\right)
\]
and, in consequence
\begin{equation}\label{split expo}
e^{\phi^{-1}(\Sig)} = \exp\left(\sum_{j=1}^n k_j\pi \left(e\p_{2j-1}^{\,2}+e\p_{2j}^{\,2}\right)\right)=\prod_{j=1}^n \exp\left[k_j\pi \left(e\p_{2j-1}^{\,2}+e\p_{2j}^{\,2}\right)\right].
\end{equation}
Let us compute $\exp\left[\pi \left(e\p_{2j-1}^{\,2}+e\p_{2j}^{\,2}\right)\right]$, $j\in \{1, \ldots, n\}$. Consider  ${\bf a}=e\p_{2j-1}-ie\p_{2j}$ and ${\bf b}=e\p_{2j-1}+ie\p_{2j}$  
 where $i$ is the usual imaginary unit in $\C$.  It is clear that ${\bf a} {\bf b}= e\p_{2j-1}^{\,2}+e\p_{2j}^{\,2}+i(e\p_{2j-1} e\p_{2j}- e\p_{2j}e\p_{2j-1}) =e\p_{2j-1}^{\,2}+e\p_{2j}^{\,2}+i$
 and $[{\bf a},{\bf b}]=2i$ which is a commuting element.  Then, 
 $\exp\left[\pi \left(e\p_{2j-1}^{\,2}+e\p_{2j}^{\,2}\right)\right]=\exp\left(\pi\, {\bf a}{\bf b}-i\pi\right)=-\exp\left(\pi\, {\bf a}{\bf b}\right)$. In order to compute $\exp\left(\pi\, {\bf a}{\bf b}\right)$ we first prove the following results.
 \begin{lem}\label{Lem xy}
For every $k\in \N$ the following relations hold.
\[ (i) \hspace{.2cm} \left[{\bf b}^k  , {\bf a} \right]=-2ik \,{\bf b}^{k-1},\hspace{1cm} (ii) \hspace{.2cm}{\bf a}^k{\bf b}^k {\bf a} {\bf b}= {\bf a}^{k+1}{\bf b}^{k+1}-2ik\, {\bf a}^k{\bf b}^k.\]
\end{lem}
\pf
\begin{itemize}
\item[$(i)$]  We proceed by induction. For $k=1$ we get $[{\bf b},{\bf a}]=-2i$ which obviously is true. Now assume that $(i)$ is true for $k\geq 1$, then for $k+1$  we get
\[{\bf b}^{k+1} {\bf a}={\bf b}\left({\bf b}^k {\bf a}\right)= {\bf b}{\bf a} {\bf b}^k-2ik \,{\bf b}^{k}=({\bf a}{\bf b}-2i){\bf b}^k-2ik \,{\bf b}^{k}={\bf a}{\bf b}^{k+1}-2i(k+1){\bf b}^{k}.\]
\item[$(ii)$] From $(i)$ we get ${\bf a}^k{\bf b}^k {\bf a} {\bf b}={\bf a}^k\left( {\bf a} {\bf b}^k-2ik \,{\bf b}^{k-1}\right) {\bf b}={\bf a}^{k+1}{\bf b}^{k+1}-2ik\, {\bf a}^k{\bf b}^k$.
$\hfill\square$
 \end{itemize}
 \begin{lem}\label{Calc}
For every $k\in \N$ it holds that $\left({\bf a} {\bf b}\right)^k=\sum_{j=1}^k (-2i)^{k-j}\, S(k,j)\, {\bf a}^j {\bf b}^j,$
where $S(n,j)$ is the Stirling number of the second kind corresponding to $k$ and $j$. 
\end{lem}
 \begin{remark}
The Stirling number of the second kind $S(k,j)$ is the number of ways of partitioning a set of $k$ elements into $j$ non empty subsets. 
We recall the following properties of the Stirling numbers,
\[S(k,1)=S(k,k)=1, \hspace{.5cm} S(k+1,j+1)=S(k,j)+(j+1)S(k,j+1),  \hspace{.5cm} \displaystyle \sum_{k=j}^\infty S(k,j) \frac{x^k}{k!}=\frac{\left(e^x-1\right)^j}{j!}. \]
\end{remark}
\noindent{\it Proof of Lemma \ref{Calc}.} \newline We proceed by induction. For $k=1$ the statement clearly is true. Now assume it to be true for $k\geq 1$. Using Lemma \ref{Lem xy}, we have  for $k+1$  that
{
\begin{align*}
({\bf a}{\bf b})^{k+1} &= \sum_{j=1}^k (-2i)^{k-j}\, S(k,j)\, {\bf a}^j {\bf b}^j {\bf a}{\bf b}  = \sum_{j=1}^k (-2i)^{k-j}\, S(k,j) {\bf a}^{j+1}{\bf b}^{j+1}+(-2i)^{k+1-j} j\, S(k,j) \, {\bf a}^j{\bf b}^j\\
&= (-2i)^k {\bf a}{\bf b}+ \left(\sum_{j=1}^{k-1}  (-2i)^{k-j} \left[S(k,j)+(j+1)S(k,j+1)\right]{\bf a}^{j+1}{\bf b}^{j+1}\right) + {\bf a}^{k+1}{\bf b}^{k+1}\\
&= \sum_{j=1}^{k+1} (-2i)^{k+1-j}\, S(k+1,j)\, {\bf a}^j {\bf b}^j,
\end{align*}
which proves the lemma.
$\hfill\square$
}

Then we obtain
\begin{align*}
e^{\pi {\bf a}{\bf b}} &= \sum_{k=0}^\infty \frac{\pi^k}{k!} ({\bf a}{\bf b})^k= 1+  \sum_{k=1}^\infty \sum_{j=1}^k \frac{\pi^k}{k!}(-2i)^{k-j}\, S(k,j)\, {\bf a}^j {\bf b}^j\\
&=1+  \sum_{j=1}^\infty  \sum_{k=j}^\infty \frac{\pi^k}{k!}(-2i)^{k-j}\, S(k,j)\, {\bf a}^j {\bf b}^j = 1+  \sum_{j=1}^\infty  (-2i)^{-j} \left[\sum_{k=j}^\infty  \frac{(-2\pi i)^k}{k!}S(k,j) \right] {\bf a}^j {\bf b}^j\\
&= 1+  \sum_{j=1}^\infty (-2i)^{-j} \frac{\left(e^{-2\pi i}-1\right)^j}{j!} {\bf a}^j {\bf b}^j\ = \ 1,
\end{align*}
from which we conclude that $\exp\left[\pi \left(e\p_{2j-1}^{\,2}+e\p_{2j}^{\,2}\right)\right]=-\exp\left(\pi\, {\bf a}{\bf b}\right)=-1$.
 \begin{remark}\label{Fou Trans}
Within the algebra $\mathcal{C}_{0,2n}=\mbox{{Alg}}_\R\{e\p_1, \ldots, e\p_{2n}\}$ the elements $e\p_{2j-1}, e\p_{2j}$ may be identified with the operators $e^{\frac{\pi}{4} i}\pa_{a_j}, e^{-\frac{\pi}{4} i}a_j$ respectively, the $a_j$'s being real variables. Indeed, these identifications immediately lead to the Weyl algebra defining relations
\begin{equation}\label{RepeprimeFourier}
e^{\frac{\pi}{4} i}\pa_{a_j}\, e^{-\frac{\pi}{4} i}a_k-e^{-\frac{\pi}{4} i}a_k\, e^{\frac{\pi}{4} i}\pa_{a_j}= \pa_{a_j} a_k-a_k\pa_{a_j}=\del_{j,k}.
\end{equation}
Hence $e\p_{2j-1}^{\,2}+e\p_{2j}^{\,2}$ may be identified with the harmonic oscillator $ i\left(\pa_{a_j}^2-a_j^2\right)$ and in consequence, the element $\exp\left[\pi \left(e\p_{2j-1}^{\,2}+e\p_{2j}^{\,2}\right)\right]$ corresponds to $\exp\left[\pi i\left(\pa_{a_j}^2-a_j^2\right)\right]$. We recall that the classical Fourier transform in one variable can be written as an operator exponential
\[\mathcal{F}[f]=\exp\left(\frac{\pi}{4} i\right) \, \exp\left(\frac{\pi}{4} i\left(\pa_{a_j}^2-a_j^2\right)\right)[f]. \]
Hence, $\exp\left[\pi i\left(\pa_{a_j}^2-a_j^2\right)\right]=-\mathcal{F}^4=-id$, where $id$ denotes the identity operator.
\end{remark}
\noindent Going back to (\ref{split expo}) we have $e^{\phi^{-1}(\Sig)} = \prod_{j=1}^n \exp\left[\pi \left(e\p_{2j-1}^{\,2}+e\p_{2j}^{\,2}\right)\right]^{k_j}=(-1)^{\sum k_j}$, 
whence $e^{B_s}=\pm 1$. Then, for $B=B_o+B_s\in \Xi_1$ such that $e^{\phi(B)}=I_{m+2n}$, we have $e^B=e^{B_o}e^{B_s}=\pm 1$, 
i.e.\  $\ker h|_{\exp(\Xi_1)}=\{-1,1\}$. In this way, we have proven the following result.
 \begin{teo}
The set $\Xi=\exp(\Xi_1) \exp(\Xi_2) \exp(\Xi_3)$ is a double covering of $\mbox{\normalfont SO}_0$.
\end{teo}
\begin{remark}\label{Fou Trans1}
As shown before, every extended superbivector of the form $B=\sum_{j=1}^n \frac{\cit_j}{2}\pi \left(e\p_{2j-1}^{\,2}+e\p_{2j}^{\,2}\right)$, $\cit_j\in\R$,
 belongs to $\Xi_1$. Then, through the identifications (\ref{RepeprimeFourier})  we can see all operators
\[\exp\left[\sum_{j=1}^n \frac{\cit_j}{2}\pi i(\pa_{a_j}^2-a_j^2)\right]=\prod_{j=1}^n \exp \left[\cit_j \frac{\pi}{2} i (\pa_{a_j}^2-a_j^2)\right]=\prod_{j=1}^n \exp\left(-\cit_j\frac{\pi}{2} i\right) \, \mathcal{F}_{a_j}^{2\cit_j},\]
as elements of the Spin group in superspace. Here, $\mathcal{F}_{a_j}^{2\cit_j}$ denotes the one-dimensional fractional Fourier transform of order $2\cit_j$ in the variable $a_j$. 
\end{remark}





\section{Conclusions}
In this paper we have shown that supervector reflections are not enough to describe the set of linear transformations leaving the inner product  invariant. This constitutes a very important difference with the classical case in which the algebra of bivectors $\underline{x}\wedge \underline{y}$ is isomorphic to the special orthogonal algebra $\mathfrak{so}(m)$. Such a property is no longer fulfilled in this setting. The real projection of the algebra of superbivectors $\R_{m|2n}^{(2)}(\Lam_N)$ does not include the symplectic algebra structure which is present in the Lie algebra of supermatrices $\mathfrak{so}_0$, corresponding to the group of super rotations.

That fact has an major impact on the definition of the Spin group in this setting. The set of elements defined through the multiplication of an even number of unit vectors in $\R^{m,2n}(\Lam_N)$ does not suffice for describing $\Spin(m|2n)(\Lam_N)$. A suitable alternative, in this case, is to define the (super) spin elements as products of exponentials of {\bf extended} superbivectors. Such an extension of the Lie algebra of superbivectors contains, through the corresponding identifications, harmonic oscillators. This way, we obtain the Spin group as a cover of the set of superrotations SO${}_0$ through the usual representation $h$.  In addition, every fractional Fourier transform can be identified with a spin element.

 In the related paper \cite{DS_Guz_Somm4} , we have already proven the invariance of the (super) Dirac operator $\pa_x$ under the corresponding actions of this Spin group.We have also studied there the invariance of the Hermitian system under the action of the corresponding Spin subgroup.
 

\section*{Acknowledgements}
Alí Guzmán Adán is supported by a BOF post doctoral grant from Ghent University.


\bibliographystyle{abbrv}

 
 

\end{document}